\documentclass[11pt,draftcls,onecolumn]{IEEEtran}
\usepackage{amsmath,amssymb,amsfonts}
\usepackage{graphics,graphicx}
\usepackage{epsfig}
\usepackage{cite}
\usepackage{bm}
\usepackage{latexsym}
\usepackage{eepic}
\usepackage{algorithmic}
\usepackage{algorithm}
\usepackage{color}
\usepackage[none]{hyphenat}
\usepackage{dsfont}

\newcommand{\mat}[1]{\ensuremath{{\mathbf{#1}}}}
\newcommand{\matc}[1]{\ensuremath{{\mathcal{#1}}}}

\newcommand{\inmat}[3][C]{\ensuremath{\in{\mathbb{#1}}^{{#2}\times{#3}}}}
\newcommand{\inmatR}[3][R]{\ensuremath{\in{\mathbb{#1}}^{{#2}\times{#3}}}}
\newcommand{\inten}[4][C]{\ensuremath{\in{\mathbb{#1}}^{{#2}\times{#3}\times{#4}}}}

\newcommand{\be}{\begin{equation}}
\newcommand{\ee}{\end{equation}}
\newcommand{\beq}{\begin{eqnarray}}
\newcommand{\eeq}{\end{eqnarray}}
\newcommand{\bi}{\begin{itemize}}
\newcommand{\ei}{\end{itemize}}
\newcommand{\khatri}{\diamond} 
\newcommand{\kron}{\otimes} 

\newcommand{\ba}{\begin{array}}
\newcommand{\ea}{\end{array}}

\newenvironment{proof} { \vspace{\smallskipamount}\par {\it Proof.}~}{\hfill $\Box$ \vspace{\medskipamount}\par }
\title{Overview of Constrained PARAFAC Models}

\author{G\'{e}rard Favier and Andr\'{e} L. F. de Almeida,

\thanks{G\'{e}rard Favier is with the I3S Laboratory, University of Nice-Sophia Antipolis (UNS), CNRS, France. Andr\'{e} L. F. de Almeida is with the Wireless Telecom Research Group, Department of Teleinformatics Engineering, Federal University of Cear\'{a}, Fortaleza, Brazil. E-mails: favier@i3s.unice.fr,\ andre@gtel.ufc.br.}}

\begin{document}

\maketitle

\begin{abstract} 
In this paper, we present an overview of constrained PARAFAC models where the constraints model linear dependencies among columns of the factor matrices of the tensor decomposition, or  alternatively, the pattern of interactions between  different modes of the tensor  which are captured by the equivalent core tensor. Some tensor prerequisites with a particular emphasis on mode combination using Kronecker products of canonical vectors that makes easier matricization operations, are first introduced. This Kronecker product based approach is also formulated in terms of the index notation, which provides an original and concise formalism for both matricizing tensors and writing tensor models. Then, after a brief reminder of PARAFAC and Tucker models, two families of constrained tensor models, the co-called PARALIND/CONFAC and PARATUCK models, are described in a unified framework, for $N^{th}$ order tensors. New tensor models, called nested Tucker models and block PARALIND/CONFAC models, are also introduced. A link between PARATUCK models and constrained PARAFAC models is then established. Finally, new uniqueness properties of PARATUCK models are deduced from sufficient conditions for essential uniqueness of their associated constrained PARAFAC models.
\end{abstract}

\begin{keywords}
Constrained PARAFAC, PARALIND/CONFAC, PARATUCK, Tensor models, Tucker models.
\end{keywords}

\section{Introduction}

Tensor calculus was introduced in differential geometry, at the end of the $19^{th}$ century, and then tensor analysis was developed in the context of Einstein's theory of general relativity, with the introduction of index notation, the so-called Einstein summation convention, at the beginning of the $20^{th}$ century, which allows to simplify and shorten physics equations involving tensors. Index notation is also useful for simplifying multivariate statistical calculations, particularly those involving cumulant tensors \cite{Cullagh}. Generally speaking, tensors are used in physics and differential geometry for characterizing the properties of a physical system, representing fundamental laws of physics, and defining geometrical objects whose components are functions. When these functions are defined over a continuum of points of a mathematical space, the tensor forms what is called a tensor field, a generalization of vector field used to solve problems involving curved surfaces or spaces, as it is the case of curved space-time in general relativity. From a mathematical point of view, two other approaches are possible for defining tensors, in terms of tensor products of vector spaces, or multilinear maps. Symmetric tensors can also be linked with homogeneous polynomials \cite{Como02:oxford}.

After first tensor developments by mathematicians and physicists, the need of analysing collections of data matrices that can be seen as three-way data arrays, gave rise to three-way models for data analysis, with the pioneering works of Tucker (1966) in psychometrics \cite{Tucker66}, and Harshman (1970) in phonetics \cite{Harshman70}, who proposed what is now referred to as the Tucker and the PARAFAC {\color{black}(parallel factor) decompositions/models, respectively. PARAFAC decompositions were independently proposed by Carroll and Chang in 1970 \cite{CarrollChang70} under the name CANDECOMP (canonical decomposition), then called CP (for CANDECOMP/PARAFAC) in \cite{Kiers2000}.} For an history of the development of multi-way models in the context of data analysis, see \cite{Kroonenberg2008}. Since the nineties, multi-way analysis has known a growing success in chemistry and especially in chemometrics. See Bro's thesis (1998) \cite{Bro98} and the book by Smilde et al. (2004) \cite{Smilde_book} for a description of various chemical applications of three-way models, with a pedagogical presentation of these models and of various algorithms for estimating their parameters. At the same period, tensor tools were developed for signal processing applications, more particularly for solving the so-called blind source separation (BSS) problem using cumulant tensors. See \cite{Cardoso90}, \cite{CardoComon90}, \cite{Cardoso91}, and De Lathauwer's thesis \cite{DeLathau97} where the concept of HOSVD (high order singular value decomposition) is introduced, a tensor tool generalizing the standard matrix SVD to arrays of order higher than two. A recent overview of BSS approaches and applications can be found in the handbook co-edited by Comon and Jutten \cite{CJ10}.

Nowadays, (high order) tensors, also called multi-way arrays in the data analysis community, play an important role in many fields of application for representing and analysing multidimensional data, as in psychometrics, chemometrics, food industry, environmental sciences, signal/image processing, computer vision, neuroscience, information sciences, data mining, pattern recognition, among many others. Then, they are simply considered as multidimensional arrays of numbers, constituting a generalization of vectors and matrices that are first- and second-order tensors respectively, to orders higher than two. Tensor models, also called tensor decompositions, are very useful for analysing multidimensional data under the form of signals, images, speech, music sequences, or texts, and also for designing new systems as it is the case of wireless communication systems since the publication of the seminal paper by Sidiropoulos et al., in 2000 \cite{Sid00}. Besides the references already cited, overviews of tensor tools, models, algorithms, and applications can be found in \cite{Cichocki09}, \cite{Kolda09}, \cite{Acar09}, \cite{Morten11}.

Tensor models incorporating constraints (sparsity; non-negativity; smoothness; symmetry; column orthonormality of factor matrices; Hankel, Toeplitz, and Vandermonde structured matrix factors; allocation constraints,...) have been the object of intensive works, during the last years. Such constraints can be inherent to the problem under study, or the result of a system design. An overview of constraints on components of tensor models most often encountered in multi-way data analysis can be found in \cite{Kroonenberg2008}. {\color{black} Incorporation of constraints in tensor models may facilitate physical interpretabibility of matrix factors. Moreover, imposing constraints may allow to relax uniqueness conditions, and to develop specialized parameter estimation algorithms with improved performance both in terms of accuracy and computational cost, as it is the case of CP models with a columnwise orthonormal factor matrix \cite{SLC2012}.}
One can classify the constraints into three main categories: i) sparsity/non-negativity, ii) structural, iii) linear dependencies/mode interactions. It is worth noting that the three categories of constraints involve specific parameter estimation algorithms, the first two ones generally inducing an improvement of uniqueness property of the tensor decomposition, while the third category implies a reduction of uniqueness, named partial uniqueness. We briefly review the main results concerning the first two types of constraints, section III of this paper being dedicated to the third category.

Sparse and non-negative tensor models have recently been the subject of many works in various fields of applications like computer vision (\cite{Shashua2005}, \cite{Hazan2005}), image compression \cite{Friedlander2008}, hyperspectral imaging \cite{Zhang2008}, music genre classification \cite{Benetos2010} and audio source separation \cite{Ozerov2011}, multi-channel EEG (electroencephalography) and network traffic analysis \cite{Acar2010}, fluorescence analysis \cite{Royer2011}, data denoising and image classification \cite{Phan2011}, among many others. Two non-negative tensor models have been more particularly studied in the literature, the so-called non-negative tensor factorization (NTF), i.e. PARAFAC models with non-negativity constraints on the matrix factors, and non-negative Tucker decomposition (NTD), i.e. Tucker models with non-negativity constraints on the core tensor and/or the matrix factors. The crucial importance of NTF/NTD for multi-way data analysis applications results from the very large volume of real-world data to be analyzed under constraints of sparseness and non-negativity of factors to be estimated, when only non-negative parameters are physically interpretable. Many NTF/NTD algorithms are now available. Most of them can be viewed as high-order extensions of non-negative matrix factorization (NMF) methods, in the sense that they are based on an alternating minimization of cost functions incorporating sparsity measures (also named distances or divergences) with application of NMF methods to matricized or vectorized forms of the tensor to be decomposed. See for instance \cite{Welling2001}, \cite{Friedlander2008}, \cite{Cichocki09}, \cite{Royer2011} for NTF, and \cite{Morup2008}, \cite{Phan2011} for NTD. An overview of NMF and NTF/NTD algorithms can be found in \cite{Cichocki09}.

The second category of constraints concerns the case where the core tensor and/or some matrix factors of the tensor model have a special structure. For instance, we recently proposed a nonlinear CDMA scheme for multiuser SIMO communication systems that is based on a constrained block-Tucker2 model whose core tensor, composed of the information symbols to be transmitted and their powers up to a certain degree, is characterized by matrix slices having a Vandermonde or a Hankel structure \cite{FB09EUSIPCO}, \cite{FBA12}. We also developed Volterra-PARAFAC models for nonlinear system modeling and identification. These models are obtained by expanding high-order Volterra kernels, viewed as symmetric tensors, by means of symmetric or doubly symmetric PARAFAC decompositions \cite{FKB2012}, \cite{BF2012}. Block structured nonlinear systems like Wiener, Hammerstein, and parallel-cascade Wiener systems, can be identified from their associated Volterra kernels that admit symmetric PARAFAC decompositions with Toeplitz factors \cite{KF2009}, \cite{Favier2009}.
Symmetric PARAFAC models with Hankel factors, and symmetric block PARAFAC models with block Hankel factors are encountered for blind identification of MIMO linear channels using fourth-order cumulant tensors, in the cases of memoryless and convolutive channels, respectively \cite{Fernandes2007}, \cite{Fernandes2009}. In the presence of structural constraints, specific estimation algorithms can be derived as it is the case for {\color{black}symmetric CP decompositions  \cite{Brachat2010}, CP decompositions with Toeplitz factors (in \cite{Nion2008} an iterative solution was proposed, whereas in \cite{KF2009b} a non-iterative algorithm was developed), Vandermonde factors \cite{SL2013}, circulant factors \cite{GoF2014}, or more generally with banded and/or structured matrix factors \cite{Comon2010}, \cite{SC2013},} and for Hankel and Vandermonde structured core tensors \cite{FBA12}.

The rest of this paper is organized as follows. In Section II, we present some tensor prerequisites with a particular emphasis on mode combination using Kronecker products of canonical vectors that makes easier the matricization operations, especially to derive matrix representations of tensor models. {\color{black}This Kronecker product based approach is also formulated in terms of the index notation, which provides an original and concise formalism for both matricizing tensors and writing tensor models.} We also present the two most common tensor models, the so called Tucker and PARAFAC models, in a general framework, i.e. for $N^{th}$-order tensors. Then, in Section III, two families of constrained tensor models, the co-called PARALIND/CONFAC and PARATUCK models, are described in a unified way, with a generalization to $N^{th}$ order tensors. New tensor models, called nested Tucker models and block PARALIND/CONFAC models, are also introduced. A link between PARATUCK models and constrained PARAFAC models is also established. In Section IV, uniqueness properties of PARATUCK models are deduced using this link. The paper is concluded in Section V.\\


\subsection*{\textit{Notations and definitions:}}

$\mathbb{R}$ and $ \mathbb{C} $ denote the fields of real and complex numbers, respectively. Scalars, column vectors, matrices, and high order tensors are denoted by lowercase, boldface lowercase, boldface uppercase, and calligraphic letters, e.g. $a$, $\mat a$, $\mat A$, and $\matc A$, respectively. The vector $\mat A_{i.}$ (resp. $\mat A_{.j}$) represents the $i^{th}$ row (resp. $j^{th}$ column) of $\mat A$.

$\mat I_{N}$, $\mat 1_{N}^{T}$, and $\mat e_{n}^{(N)}$ stand for the identity matrix of order $N$, the all-ones row vector of dimensions ${1 \times N}$, and the $n^{th}$ canonical vector of the Euclidean space $\mathbb{R}^{N}$, respectively.

$\mat A^{T}$, $\mat A^{H}$, $\mat A^{\dag}$, tr$(\mat A)$, and $r_{\mat A}$ denote the transpose, the conjugate (Hermitian) transpose, the Moore-Penrose pseudo-inverse, the trace, and the rank of $\mat A$, respectively. $D_{i}(\mat A)=diag(\mat A_{i.})$ represents the diagonal matrix having the elements of the $i^{th}$ row of $\mat A$ on its diagonal. The operator $bdiag(.)$ forms a block-diagonal matrix from its matrix arguments, while the operator vec(.) transforms a matrix into a column vector by stacking the columns of its matrix argument one on top of the other one. In case of a tensor $\matc X$, the vec operation is defined in (\ref{vectorization}).
\

The outer product (also called tensor product), and the matrix Kronecker, Khatri-Rao (column-wise Kronecker), and Hadamard (element-wise) products are denoted by $\circ$, $\otimes$, $\diamond$, and $\odot$, respectively.
\

Let us consider the set $\mathds{S}=\{n_1,\ldots,n_N \}$ obtained by permuting the elements of the set $\{1,\ldots,N\}$. For $\mat A^{(n)}\inmat{I_{n}}{R_{n}}$ and $\mat u^{(n)}\inmat{I_{n}}{1}$, $n=1,\cdots ,N$, we define
\beq
&&\hspace{-3ex}\underset{n \in \mathds{S}}{\kron}\mat A^{(n)}=\mat A^{(n_1)}\kron \mat A^{(n_2)}\kron \cdots \kron \mat A^{(n_N)}\inmat{I_{n_{1}}\cdots I_{n_{N}}}{R_{n_{1}}\cdots R_{n_{N}}}; \label{defkron}\\
&&\hspace{-3ex}\underset{n \in \mathds{S}}{\khatri}\mat A^{(n)}=\mat A^{(n_1)}\khatri \mat A^{(n_2)}\khatri \cdots \khatri \mat A^{(n_N)}\inmat{I_{n_{1}}\cdots I_{n_{N}}}{R}, \nonumber\\
&& \quad \textrm{when} \quad R_n=R, \forall n=1,\cdots ,N ;\label{defkhatri}\\
&&\hspace{-3ex}\underset{n \in \mathds{S}}{\odot}\mat A^{(n)}=\mat A^{(n_1)}\odot \mat A^{(n_2)}\odot \cdots \odot \mat A^{(n_N)}\inmat{I}{R},\nonumber\\
&&\quad \textrm{when} \, I_n=I, \, \textrm{and} \, R_n=R, \, \forall n=1,\cdots ,N; \nonumber\\
&&\hspace{-3ex}\underset{n \in \mathds{S}}{\circ}\mat u^{(n)}=\mat u^{(n_1)}\circ \mat u^{(n_2)}\circ \cdots \circ \mat u^{(n_N)} \inten{I_{n_{1}}}{\cdots}{I_{n_{N}}}. \nonumber
\eeq
The outer product of $N$ non-zero vectors defines a rank-one tensor of order $N$.
\

By convention, the order of dimensions is directly related to the order of variation of the associated indices. For instance, in (\ref{defkron}) and (\ref{defkhatri}), the product $I_{n_{1}}I_{n_{2}}\cdots I_{n_{N}}$ of dimensions means that $n_1$ is the index varying the most slowly while $n_N$ is the index varying the most fastly in the Kronecker products computation.

For $\mathds{S}=\{1,\ldots,N\}$, we have the following identities
\beq \label{eq: outer product def}
&&\hspace{-8ex}{\left( \underset{n \in \mathds{S}}{\circ} \mat u^{(n)} \right)}_{i_{1},\cdots, i_{N}}={\left( \overset{N}{\underset{n=1}{\circ}}\mat u^{(n)} \right)}_{i_{1},\cdots, i_{N}}=\prod\limits_{n=1}^{N}u^{(n)}_{i_{n}}, \nonumber\\
\label{eq: Kron product def}
&&\hspace{-8ex}{\left( \underset{n \in \mathds{S}}{\otimes} \mat u^{(n)} \right)}_i={\left( \overset{N}{\underset{n=1}{\otimes}} \mat u^{(n)} \right)}_i=\prod\limits_{n=1}^{N}u^{(n)}_{i_{n}} \quad \textrm{with} \quad i=i_N+\sum\limits_{n=1}^{N-1}(i_n-1)\prod\limits_{j=n+1}^{N}I_j.
\eeq

In particular, for {\mat u\inmat{I}{1}}, {\mat v\inmat{J}{1}}, {\mat w\inmat{K}{1}}
\beq \label{eq: outer product prop}
&&{\matc X=\mat u\circ \mat v\circ \mat w} \inten{I}{J}{K}\Leftrightarrow x_{ijk}=u_i v_j w_k,\nonumber\\
&&{\mat x=\mat u\otimes \mat v\otimes \mat w} \inmat{IJK}{1}\Leftrightarrow x_{k+(j-1)K+(i-1)JK}=u_i v_j w_k.\nonumber
\eeq


\
Some useful matrix formulae are recalled in the Appendix.

\section{Tensor Prerequisites}

In this paper, a tensor is simply viewed as a multidimensional array of measurements. Depending that these measurements are real- or complex-valued, we have a real- or complex-valued tensor, respectively. The order $N$ of a tensor refers to the number of indices that characterize its elements $x_{i_1,\cdots,i_N }$, each index $i_n$ ($i_n=1,\cdots,I_N, \textrm{for}\,\, n=1,\cdots, N$) being associated with a dimension, also called a way, or a mode, and $I_n$ denoting the mode-$n$ dimension.

An $N^{th}$-order complex-valued tensor $\matc X \inten{I_{1}}{\cdots}{I_{N}}$, also called an $N$-way array, of dimensions $I_{1}\times \cdots \times I_{N}$, can be written as
\beq \label{eq: tensor def}
&&{\matc X=\sum\limits_{i_1=1}^{I_1}\cdots \sum\limits_{i_N=1}^{I_N}x_{i_1,\cdots,i_N }\overset{N}{\underset{n=1}{\circ}}\mat e^{(I_n)}_{i_n}}.
\eeq
\
The coefficients $x_{i_1,\cdots,i_N }$ represent the coordinates of $\matc X$ in the canonical basis $\{\overset{N}{\underset{n=1}{\circ}}\mat e^{(I_n)}_{i_n}, i_n=1,\cdots, I_n;n=1,\cdots, N\}$ of the space $\mathbb{C}^{{I_{1}}\times{\cdots}\times{I_{N}}}$.\

The identity tensor of order $N$ and dimensions $I\times \cdots \times I$, denoted by $\matc I_{N,I}$ or simply $\matc I$, is a diagonal hypercubic tensor whose elements $\delta_{i_1,\cdots,i_N }$ are defined by means of the generalized Kronecker delta, i.e. $\delta_{i_1,\cdots,i_N }=\left\{\begin{array}{ll} 1 \quad \textrm{if}\quad i_1=\cdots=i_N\\
0 \quad \textrm{otherwise} \end{array}\right.$, and $I_n=I,\forall n=1,\cdots ,N$. It can be written as
\beq
&&\matc I_{N,I}=\sum\limits_{i=1}^{I} \underset{N\,\, \textrm{terms}}{\underbrace{\mat e^{(I)}_i \circ \cdots \circ \mat e^{(I)}_i}}.\nonumber
\eeq

Different reduced order tensors can be obtained by slicing the tensor $\matc X \inten{I_{1}}{\cdots}{I_{N}}$ along one mode or $p$ modes, i.e. by fixing one index $i_n$ or a set of $p$ indices $\{i_{n_1},\ldots,i_{n_p} \}$, which gives a tensor of order $N-1$ or $N-p$, respectively. For instance, by slicing $\matc X$ along its mode-$n$, we get the $i^{th}_n$ mode-$n$ slice of $\matc X$, denoted by $\matc X_{\ldots i_n \ldots}$, that can be written as

\beq
{\matc X_{\ldots i_n \ldots}=\sum\limits_{i_1=1}^{I_1} \cdots \sum\limits_{i_{n-1}=1}^{I_{n-1}} \sum\limits_{i_{n+1}=1}^{I_{n+1}} \cdots \sum\limits_{i_N=1}^{I_N} x_{i_{1},\cdots ,i_{n}, \cdots,i_{N}} \mat e^{(I_{n+1})}_{i_{n+1}} \circ \cdots \circ \mat e^{(I_N)}_{i_N} \circ \mat e^{(I_1)}_{i_1} \circ \cdots}\nonumber \\
 {\cdots \circ \mat e^{(I_{n-1})}_{i_{n-1}}\inten{I_{n+1}}{\cdots \times I_N \times I_1\times \cdots} {I_{n-1}}}.\nonumber
\eeq
\
For instance, by slicing the third-order tensor $\matc X \inten{I}{J}{K}$ along each mode, we get three types of matrix slices, respectively called horizontal, lateral, and frontal slices:
\beq
\mat X_{i..}\inmat{J}{K},  \mat X_{.j.}\inmat{K}{I} \,\, \textrm{and} \,\, \mat X_{..k}\inmat{I}{J}, \nonumber\\
\textrm{with} \,\,  i=1,\ldots,I;  j=1,\ldots,J;  k=1,\ldots,K.\nonumber
\eeq
{\color{black}
\subsection{\bf{Tensor Hadamard Product}}

{Consider $\matc A\inten{R_{1}\times \cdots}{ R_N}{I_{1}\times \cdots \times I_{P_1}}$ and $\matc B\inten{R_{1}\times \cdots}{ R_N}{I_{P_1+1}\times \cdots \times I_P}$, and the ordered subset $\mathds{R}=\{r_1,\cdots,r_N\}$.
The Hadamard product of $\matc A$ with $\matc B$ along their common modes, gives a tensor $\matc C\inten{R_{1}\times \cdots}{ R_N}{I_1\times \cdots \times I_P}$ such that}
\beq
\matc C=\,\matc A\ {\underset{\mathds{R}}{\odot}}\, \matc B \Leftrightarrow c_{r_1,\cdots,r_N,i_1,\cdots,i_P}=a_{r_1,\cdots,r_N,i_1,\cdots,i_{P_1}}b_{r_1,\cdots,r_N,i_{P_1+1},\cdots,i_P}\nonumber
\eeq
For instance, given two third-order tensors $\matc A\inten{R_1}{R_2}{I_1}$ and $\matc B\inten{R_1}{R_2}{I_2}$, the Hadamard product $\matc A\ {\underset{\{r_1,r_2\}}{\odot}}\, \matc B$ gives a fourth-order tensor $\matc C\inten{R_1}{R_2\times I_1}{I_2}$ such that
\beq
c_{r_1,r_2,i_1,i_2}=a_{r_1,r_2,i_1}b_{r_1,r_2,i_2}\nonumber.
\eeq
Such a tensor Hadamard product can be calculated by means of the matrix Hadamard product of extended tensor unfoldings as defined in Eq. (\ref{tensor extension 1}) and (\ref{tensor extension 2}) (see also Eq. (\ref{eq: extension formula})-(\ref{eq: extension formula c}) in the Appendix A.5). For the example above, we have
\beq
\mat C_{R_1R_2\times I_1I_2}=\mat A_{R_1R_2\times I_1}(\mat I_{I_1}\otimes \mat 1_{I_2}^T)\odot\mat B_{R_1R_2\times I_2}(\mat 1_{I_1}^T\otimes \mat I_{I_2})\nonumber
\eeq
\subsection*{Example:} For $\mat A_{R\times I_1}=\left[\begin{array}{cc} a_1 & a_2 \\
a_3 & a_4 \end{array}\right]$, $\mat B_{R\times I_2}=\left[\begin{array}{cc} b_1 & b_2 \\
b_3 & b_4 \end{array}\right]$, and the tensor $\matc C$ such as $c_{r,i_1,i_2}=a_{r,i_1}b_{r,i_2}$, a mode-1 flat matrix unfolding of $\matc C$ is given by
\beq
\mat C_{R\times I_1I_2}&=&\mat A_{R\times I_1}(\mat I_2\otimes \mat 1_2^T)\odot\mat B_{R\times I_2}(\mat 1_2^T\otimes \mat I_2)\nonumber\\
&=&\left[\begin{array}{cccc} a_1 & a_1 & a_2 & a_2  \\
a_3 & a_3 & a_4 & a_4 \end{array}\right]\odot \left[\begin{array}{cccc} b_1 & b_2 & b_1 & b_2  \\
b_3 & b_4 & b_3 & b_4 \end{array}\right]\nonumber\\
&=&\left[\begin{array}{cccc} a_1 b_1 & a_1b_2 & a_2b_1 & a_2b_2  \\
a_3b_3 & a_3b_4 & a_4b_3 & a_4b_4 \end{array}\right]\nonumber
\eeq
}
\subsection{\bf{Mode Combination}}

Different contraction operations can be defined depending on the way according to which the modes are combined. Let us partition the set $\{1,\ldots,N\}$ in $N_1$ ordered subsets $\mathds{S}_{n_1}$, constituted of $p(n_1)$ elements with $\sum\limits_{n_1=1}^{N_1}p(n_1)=N$. Each subset $\mathds{S}_{n_1}$ is associated with a combined mode of dimension $J_{n_1}=\underset{n \in \mathds{S}_{n_1}}{\prod I_n}$. These mode combinations allow to rewrite the $N^{th}$-order tensor $\matc X \inten{I_{1}}{\cdots}{I_{N}}$ under the form of an $N^{th}_1$-order tensor $\matc Y \inten{J_{1}}{\cdots}{J_{N_1}}$ as follows
\beq \label{eq: contraction def}
\matc Y=\sum\limits_{j_1=1}^{J_1}\cdots \sum\limits_{j_{N_1}=1}^{J_{N_1}}x_{j_1,\cdots,j_{N_1} }\overset{N_1}{\underset{n_1=1}{\circ}}\mat e^{(J_{n_1})}_{j_{n_1}} \, \, \textrm{with}\, \, \mat e^{(J_{n_1})}_{j_{n_1}}=\underset{n \in \mathds{S}_{n_1}}{\otimes} \mat e^{(I_n)}_{i_n}.
\eeq

Two particular mode combinations corresponding to the vectorization and matricization operations are now detailed.

\subsection{\bf{Vectorization}}

The vectorization of $\matc X \inten{I_{1}}{\cdots}{I_{N}}$ is associated with a combination of the $N$ modes into a unique mode of dimension $J=\prod\limits_{n=1}^{N}I_n$, which amounts to replace the outer product in (\ref {eq: tensor def}) by the Kronecker product
\beq \label{vectorization}
&&{\textrm{vec}(\matc X)=\sum\limits_{i_1=1}^{I_1}\cdots \sum\limits_{i_N=1}^{I_N}x_{i_1,\cdots,i_N }\overset{N}{\underset{n=1}{\otimes}}\mat e^{(I_n)}_{i_n}\inmat{I_{1} \cdots I_{N}}{1}}
\eeq
the element $x_{i_1,\cdots,i_N }$ of $\matc X$ being the $i^{th}$ entry of $\textrm{vec}(\matc X)$ with $i$ defined as in (\ref{eq: Kron product def}).

The vectorization can also be carried out after a permutation of indices $\pi(i_n), n=1,\cdots,N$.



\subsection{\bf{Matricization or Unfolding}}

There are different ways of matricizing the tensor $\matc X$ according to the partitioning of the set $\{1,\ldots,N\}$ into two ordered subsets $\mathds{S}_1$ and $\mathds{S}_2$, constituted of $p$ and $N-p$ indices, respectively. A general formula for the matricization is, for $p\in[1,N-1]$
\beq \label {eq: matricization formula}
&&{\mat X_{\mathds{S}_1;\mathds{S}_2}=\sum\limits_{i_1=1}^{I_1}\cdots \sum\limits_{i_N=1}^{I_N}x_{i_1,\cdots,i_N }{\left( \underset{n \in \mathds{S}_1}{\otimes}  \mat e^{(I_n)}_{i_n} \right)}{\left( \underset{n \in \mathds{S}_2}{\otimes}  \mat e^{(I_n)}_{i_n} \right)}^T \inmat{J_1}{J_2}}
\eeq
with $J_{n_1}=\underset{n \in \mathds{S}_{n_1}}{\prod I_n}$, for $n_1=1 \, \textrm{and}\,\,  2$. From (\ref{eq: matricization formula}), we can deduce the following expression of the element $x_{i_1,\cdots,i_N }$ in terms of the matrix unfolding $\mat X_{\mathds{S}_1;\mathds{S}_2}$
\beq \label {eq: tensor element}
&&{x_{i_1,\cdots,i_N }={\left( \underset{n \in \mathds{S}_1}{\otimes}  \mat e^{(I_n)}_{i_n} \right)}^T{\mat X_{\mathds{S}_1;\mathds{S}_2}}{\left( \underset{n \in \mathds{S}_2}{\otimes}  \mat e^{(I_n)}_{i_n} \right)}}.
\eeq

%

\subsection{\bf{Particular case: mode-$n$ matrix unfoldings $\mat X_n$}}

A flat mode-$n$ matrix unfolding of the tensor $\matc X$ corresponds to an unfolding of the form $\mat X_{\mathds{S}_1;\mathds{S}_2}$ with $\mathds{S}_1=\{n\}$ and $\mathds{S}_2=\{n+1,\cdots,N,1,\cdots,n-1\}$, which gives
\beq \label {eq:mode-n unfolding}
&&\mat X_{I_n\times I_{n+1}\cdots I_NI_1\cdots I_{n-1}}=\mat X_n\nonumber\\
&&=\sum\limits_{i_1=1}^{I_1}\cdots \sum\limits_{i_N=1}^{I_N}x_{i_1,\cdots,i_N }{\mat e^{(I_n)}_{i_n}}{\left( \underset{n \in \mathds{S}_2}{\otimes}  \mat e^{(I_n)}_{i_n} \right)}^T\inmat{I_n}{I_{n+1}\cdots I_NI_1\cdots I_{n-1}}.
\eeq
We can also define a tall mode-$n$ matrix unfolding of $\matc X$, by choosing $\mathds{S}_1=\{n+1,\cdots,N,1,\cdots,n-1\}$ and $\mathds{S}_2=\{n\}$. Then, we have $\mat X_{I_{n+1}\cdots I_NI_1\cdots I_{n-1}\times I_n}=\mat X_n^T\inmat{I_{n+1}\cdots I_NI_1\cdots I_{n-1}}{I_n}$.

The column vectors of a flat mode-$n$ matrix unfolding $\mat X_n$ are the mode-$n$ vectors of $\matc X$, and the rank of $\mat X_n$, i.e. the dimension of the mode-$n$ linear space spanned by the mode-$n$ vectors, is called mode-$n$ rank of $\matc X$, denoted by $\textrm{rank}_n(\matc X)$.

In the case of a third-order tensor $\matc X \inten{I}{J}{K}$, there are six different flat unfoldings, denoted $\mat X_{I\times JK}$, $\mat X_{I\times KJ}$, $\mat X_{J\times KI}$, $\mat X_{J\times IK}$, $\mat X_{K\times IJ}$, $\mat X_{K\times JI}$. For instance, we have
\beq \label{eq:flat unfolding}
\mat X_{I\times JK}=\mat X_{\{1\};\{2,3\}}=\sum\limits_{i=1}^{I}\sum\limits_{j=1}^{J}\sum\limits_{k=1}^{K}x_{i,j,k}\,\mat e^{(I)}_i(\mat e^{(J)}_j\otimes \mat e^{(K)}_k)^T.
\eeq
Using the properties (\ref {eq: kron prop1}), (\ref {eq: kron assoc prop}), and (\ref {eq: kron distribut prop}) of the Kronecker product gives
\beq
\mat X_{I\times JK}&=&\sum\limits_{j=1}^{J}(\mat e^{(J)}_j)^T\otimes\sum\limits_{i=1}^{I}\sum\limits_{k=1}^{K}x_{i,j,k}\mat e^{(I)}_i(\mat e^{(K)}_k)^T\nonumber\\
&=&\sum\limits_{j=1}^{J}(\mat e^{(J)}_j)^T\otimes(\mat X_{.j.})^T=\left[\begin{array}{ccc}\mat X^T_{.1.} & \cdots & \mat X^T_{.J.} \nonumber \end{array}\right]\inmat {I}{JK}.
\eeq
Similarly, there are six tall matrix unfoldings, denoted $\mat X_{JK\times I}$, $\mat X_{KJ\times I}$, $\mat X_{KI\times J}$, $\mat X_{IK\times J}$, $\mat X_{IJ\times K}$, $\mat X_{JI\times K}$, like for instance
\beq \label{eq:flat unfolding b}
\mat X_{JK\times I}=\sum\limits_{i=1}^{I}\sum\limits_{j=1}^{J}\sum\limits_{k=1}^{K}x_{i,j,k}\,(\mat e^{(J)}_j\otimes \mat e^{(K)}_k){\mat e^{(I)}_i}^T=\mat X^T_{I\times JK}\inmat {JK}{I}.
\eeq
Applying (\ref {eq: tensor element}) to (\ref{eq:flat unfolding}) gives
\beq
x_{i,j,k}=(\mat e^{(I)}_i)^T\mat X_{I\times JK}(\mat e^{(J)}_j\otimes \mat e^{(K)}_k)=[\mat X_{I\times JK}]_{i,(j-1)K+k}.\nonumber
\eeq

\subsection{\bf{Mode-$n$ product of a tensor with a matrix or a vector}}

The mode-$n$ product of $\matc X \inten{I_{1}}{\cdots}{I_{N}}$ with $\mat A\inmat{J_n}{I_n}$ along the $n^{th}$ mode, denoted by $\matc X {\times}_n \mat A$, gives the tensor $\matc Y$ of order $N$ and dimensions $I_{1}\times \cdots \times I_{n-1} \times  J_{n}\times I_{n+1} \times \cdots \times I_N$, such as \cite{CPK80}
\beq \label{eq: mode-n prod def}
&&{y_{i_1,\cdots,i_{n-1} ,j_n ,i_{n+1} , \cdots, i_N }=\sum\limits_{i_n=1}^{I_n}a_{{j_n},{i_n}} x_{i_1,\cdots,i_{n-1} ,i_n ,i_{n+1} , \cdots, i_N }}
\eeq
which can be expressed in terms of mode-$n$ matrix unfoldings of $\matc X$ and $\matc Y$
\beq
&&{\mat Y_n=\mat A \mat X_n}.\nonumber
\eeq
This operation can be interpreted as the linear map from the mode-$n$ space of $\matc X$ to the mode-$n$ space of $\matc Y$, associated with the matrix $\mat A$.
\

The mode-$n$ product of $\matc X \inten{I_{1}}{\cdots}{I_{N}}$ with the row vector $\mat u^T \inmat{1}{I_n}$ along the $n^{th}$ mode, denoted by $\matc X {\times}_n \mat u^T$, gives a tensor $\matc Y$ of order $N-1$ and dimensions $I_{1}\times \cdots \times I_{n-1} \times I_{n+1} \times \cdots \times I_N$, such as
\beq
&&{y_{i_1,\cdots,i_{n-1} ,i_{n+1} , \cdots, i_N }=\sum\limits_{i_n=1}^{I_n}{u_{i_n}} x_{i_1,\cdots,i_{n-1} ,i_n ,i_{n+1} , \cdots, i_N }}\nonumber
\eeq
\
that can be written in vectorized form as ${\textrm{vec}}^T(\matc Y)=\mat u^T\mat X_n\inmat{1}{I_{n+1}\cdots I_NI_1\cdots I_{n-1}}$.
\

When multiplying a $N^{th}$-order tensor by row vectors along $p$ different modes, we get a tensor of order $N-p$. For instance, for a third-order tensor $\matc X \inten{I}{J}{K}$, we have
\beq
\mat x_{ij.}=\matc X {\times}_1 \,{\mat e^{(I)}_{i}}^T{\times}_2 \,{\mat e^{(J)}_{j}}^T, \quad x_{ijk}=\matc X {\times}_1 \,{\mat e^{(I)}_{i}}^T {\times}_2  \,{\mat e^{(J)}_{j}}^T {\times}_3  \,{\mat e^{(K)}_{k}}^T.\nonumber
\eeq


Considering an ordered subset $\mathds{S}=\{m_1,\ldots,m_P \}$ of the set $\{1,\ldots,N\}$, a series of mode-$m_p$ products of $\matc X \inten{I_{1}}{\cdots}{I_{N}}$ with $\mat A^{(m_p)}\inmat{J_{m_p}}{I_{m_p}}$, $p\in\{1,\ldots,P\}$, $P\leq N$, will be concisely noted as
\beq
&&{\matc X {\times}_{m_1} {\mat A^{(m_1)}} \cdots {\times}_{m_P}{\mat A^{(m_P)}}=\matc X  {{\times}^{m_P}_{m=m_1}} {\mat A^{(m)}}}.\nonumber
\eeq

\subsection*{Properties}

\begin{itemize}
\item For any permutation $\pi(.)$ of $P$ distinct indices $m_p\in\{1,\cdots,N\}$ such as $q_p=\pi(m_p)$, $p\in\{1,\ldots,P\}$, with $P\leq N$, we have
\beq
&&{\matc X  {{\times}^{q_{P}}_{q=q_1}} {\mat A^{(q)}}=\matc X  {{\times}^{m_{P}}_{m=m_1}} {\mat A^{(m)}}}\nonumber
\eeq
which means that the order of the mode-$m_p$ products is irrelevant when the indices $m_p$ are all distinct.
\end{itemize}
\begin{itemize}
\item For two products of $\matc X \inten{I_{1}}{\cdots}{I_{N}}$ along the same mode-$n$, with $\mat A\inmat{J_n}{I_n}$ and $\mat B\inmat{K_n}{J_n}$, we have \cite{DeLathau97}
\beq \label{eq: mode-n prod property}
&&{\matc Y=\matc X  {\times}_n {\mat A} {\times}_n {\mat B}=\matc X  {\times}_n{({\mat B}{\mat A})}\inten{I_1\times\cdots\times{I_{n-1}}}{K_n}{I_{n+1}\times\cdots\times{I_N}}}.
\eeq
\end{itemize}

{\color{black}

\subsection{Kronecker products based approach using index notation}

In this subsection, we propose to reformulate our Kronecker products based approach for tensor matricization  in terms of the index notation introduced in \cite{P2011}. Using this index notation, a column vector $\mat u\inmat{I}{1}$, a row vector $\mat v^T\inmat{1}{J}$, and a matrix $\mat X\inmat{I}{J}$ are respectively written as follows
\beq
\mat u&=&\sum\limits_{i=1}^{I}u_i\mat e_i^{(I)}=u_i\mat e_i\nonumber\\
\mat v^T&=&\sum\limits_{j=1}^{J}v_j(\mat e_j^{(J)})^T=v_j\mat e^j\nonumber\\
\mat X&=&\sum\limits_{i=1}^{I}\sum\limits_{j=1}^{j}x_{ij}(\mat e_i^{(I)}\otimes (\mat e_j^{(J)})^T)=x_{ij}\mat e_i^j\nonumber
\eeq
As with Einstein summation convention, the index notation allows to drop summation signs. If an index $i\in[1,I]$ is repeated in an expression (or more generally in a term of an equation), it means that this expression (or this term) must be summed over that index from 1 to $I$.\\

Using the index notation, the horizontal, lateral, and frontal slices of a third-order tensor $\matc X \inten{I}{J}{K}$ can be written as
\beq
\mat X_{i..}=x_{ijk}\mat e_j^k \,\,;\,\, \mat X_{.j.}=x_{ijk}\mat e_k^i \,\,;\,\, \mat X_{..k}=x_{ijk}\mat e_i^j.\nonumber
\eeq

The Kronecker products $\mat u \otimes \mat v$ and $\mat A \otimes \mat B$, with $\mat A\inmat{I}{J}$ and $\mat B\inmat{K}{L}$, can be concisely written as
\beq
\mat u \otimes \mat v&=&(u_i\mat e_i)\otimes (v_j\mat e_j)= u_iv_j\mat e_{ij}\nonumber\\
\mat A \otimes \mat B&=&(a_{ij}\mat e_i^j)\otimes (b_{kl}\mat e_k^l)=a_{ij}b_{kl}\mat e_{ik}^{jl}\nonumber
\eeq
We have also
\beq
\mat u \otimes \mat v^T&=&u_iv_j\mat e_{i}^j\nonumber\\
\mat u^T \otimes \mat v^T&=&u_iv_j\mat e^{ij}\nonumber\\
\mat A^T \otimes \mat B^T&=&a_{ji}b_{lk}\mat e_{jl}^{ik}\nonumber
\eeq
and for $\mat U=[\mat u^{(1)}\cdots \mat u^{(N)}]\inmat{I}{R}$ and $\mat V=[\mat v_1\cdots \mat v_N]\inmat{J}{R}$
\beq
\mat U \mat V^T =\mat u^{(n)} (\mat v^{(n)})^T = u_i^{(n)}v_j^{(n)}\mat e_i^j\label{UV ind not}
\eeq
Using this formalism, the Khatri-Rao product $\mat A \diamond \mat B$ can be written as follows
\beq
\mat A \diamond \mat B&=&a_{ik}b_{jk}\mat e_{ij}^{k}\nonumber\\
(\mat A \diamond \mat B)^T&=&a_{ik}b_{jk}\mat e_{k}^{ij}\label{Khatri index notation}
\eeq
Considering the set $\mathds{S}=\{n_1,\ldots,n_N \}$ obtained by permuting the elements of $\{1,\ldots,N\}$, and noting $\mat e_\mathds{I}$ the Kronecker product $\underset{n \in \mathds{S}}{\kron}\mat e_{i_n}^{(I_n)}$, with $\mathds{I}=\{i_{n_1},\cdots,i_{n_N}\}$, we have
\beq \label{veckhatri ind not}
\underset{n \in \mathds{S}}{\diamond}\mat u^{(n)}=\underset{n \in \mathds{S}}{\prod}u_{i_n}^{(n)}\mat e_{\mathds{I}}
\eeq
The Kronecker and Khatri-Rao products defined in (\ref{defkron}) and (\ref{defkhatri}), with $a_{i_n,r_n}^{(n)}$ as entry of $\mat A^{(n)}$, can then be defined as
\beq \label{Kron ind not}
\underset{n \in \mathds{S}}{\kron}\mat A^{(n)}&=&\underset{n \in \mathds{S}}{\prod}a_{i_n,r_n}^{(n)}\mat e_{i_{n_1},\cdots,i_{n_N}}^{r_{n_1},\cdots,r_{n_N}}=\underset{n \in \mathds{S}}{\prod}a_{i_n,r_n}^{(n)}\mat e_{\mathds{I}}^{\mathds{R}}\\
\underset{n \in \mathds{S}}{\diamond}\mat A^{(n)}&=&\underset{n \in \mathds{S}}{\prod}a_{i_n,r}^{(n)}\mat e_{i_{n_1},\cdots,i_{n_N}}^{r}=\underset{ n \in \mathds{S}}{\prod}a_{i_n,r}^{(n)}\mat e_{\mathds{I}}^{r}
\eeq
where $\mathds{R}=\{r_{n_1},\cdots,r_{n_N}\}$.\\
Applying these results, the unfoldings (\ref {eq: matricization formula}), (\ref{eq:flat unfolding}) and (\ref{eq:flat unfolding b}), and the formula (\ref {eq: tensor element}) can be rewritten respectively as
\beq
\mat X_{\mathds{S}_1;\mathds{S}_2}&=&x_{i_1,\cdots ,i_N}\mat e_{\mathds{I}_1}^{\mathds{I}_2}\label{unfol ind not}\\
\mat X_{I\times JK}&=&x_{i,j,k}\mat e_i^{jk}\nonumber\\
\mat X_{JK\times I}&=&x_{i,j,k}\mat e_{jk}^{i}\nonumber\\
x_{i_1,\cdots ,i_N}&=&e^{\mathds{I}_1}\mat X_{\mathds{S}_1;\mathds{S}_2}\mat e_{\mathds{I}_2}\label{unfolinverse ind not}
\eeq
where $\mathds{I}_1$ and $\mathds{I}_2$ represent the sets of indices $i_n$ associated with the sets $\mathds{S}_1$ and $\mathds{S}_2$ of index $n$, respectively.\\
We can also use the index notation for deriving matrix unfoldings of tensor extensions of a matrix $\mat B\inmat{I}{J}$. For instance, if we define the tensor $\matc A\inten{I}{J}{K}$ such as $a_{i,j,k}=b_{i,j}$ for $k=1,\cdots,K$, mode-1 flat unfoldings of $\matc A$ are given by
\beq
\mat A_{I\times JK}&=&a_{i,j,k}\mat e_i^{jk}=b_{i,j}\mat e_i^j\otimes \sum\limits_{k=1}^K\mat e^k\nonumber\\
&=&\mat B\otimes \mat 1_K^T=\mat B(\mat I_J \otimes \mat 1_K^T)\label{tensor extension 1}\\
\mat A_{I\times KJ}&=&a_{i,j,k}\mat e_i^{kj}=\sum\limits_{k=1}^K\mat e^k\otimes b_{i,j}\mat e_i^j \nonumber\\
&=&\mat 1_K^T\otimes \mat B =\mat B(\mat 1_K^T\otimes \mat I_J )\label{tensor extension 2}
\eeq
These two formulae will be used later for establishing the link between PARATUCK-(2,4) models and constrained PARAFAC-4 models. See the Appendix A.4. It is worth noting two differences between the index notation used in this paper and Einstein summation convention: $(i)$ each index can be repeated more than twice in any expression; $(ii)$ the index notation can be used with ordered sets of indices.

\subsection{Basic Tensor Models}}
We now present the two most common tensor models, i.e. the Tucker \cite{Tucker66} and PARAFAC \cite{Harshman70} models.
In \cite{Kroonenberg2008}, these models are introduced in a constructive way, in the context of three-way data analysis. The Tucker models are presented as extensions of the matrix singular value decomposition (SVD) to three-way arrays, which gave rise to the generalization as HOSVD (\cite{DeLathau97},\cite{DeLathau00}), whereas the PARAFAC model is introduced by emphasizing the Cattell's principle of parallel proportional profiles \cite{Cattell44} that underlies this model, so explaining the acronym PARAFAC. In the following, we adopt a more general presentation for multi-way arrays, i.e. tensors of any order $N$.\\

\subsubsection{\bf{Tucker Models}}

For a $N^{th}$-order tensor $\matc X \inten{I_{1}}{\cdots}{I_{N}}$, a Tucker model is defined in an element-wise form as
\beq \label{eq: Tucker model a}
&&{x_{i_1,\cdots,i_N }=\sum\limits_{r_1=1}^{R_1}\cdots \sum\limits_{r_N=1}^{R_N}g_{r_1,\cdots,r_N }\prod\limits_{n=1}^{N}a^{(n)}_{i_n,r_n}}
\eeq
with $i_n=1,\cdots,I_n$ for $n=1,\cdots,N$, where $g_{r_1,\cdots,r_N }$ is an element of the core tensor $\matc G \inten{R_{1}}{\cdots}{R_{N}}$ and $a^{(n)}_{i_n,r_n}$ is an element of the matrix factor $\mat A^{(n)}\inmat{I_n}{R_n}$.
\
{\color{black}
Using the index notation, and defining the set of indices $\mathds{R}=\{r_{n_1},\cdots,r_{n_N}\}$, the Tucker model can also be written simply as
\beq \label{Tucker ind not}
x_{i_1,\cdots,i_N}=g_{r_1,\cdots,r_N }\underset{\mathds{R}}{\prod}a^{(n)}_{i_n,r_n}
\eeq
}
Taking the definition (\ref{eq: tensor def}) into account, and noting that $\sum\limits_{i_n=1}^{I_n}a^{(n)}_{i_n,r_n}\mat e^{(I_n)}_{i_n}=\mat A^{(n)}_{.r_n}$, this model can be written as a weighted sum of $\prod\limits_{n=1}^{N}R_n$ outer products, i.e. rank-one tensors
\beq \label{eq: Tucker model b}
\matc X&=&\sum\limits_{r_1=1}^{R_1}\cdots \sum\limits_{r_N=1}^{R_N}g_{r_1,\cdots,r_N }\overset{N}{\underset{n=1}{\circ}}\mat A^{(n)}_{.r_n}\nonumber\\
&=&{\color{black}g_{r_1,\cdots,r_N }{\underset{\mathds{R}}{\circ}}\mat A^{(n)}_{.r_n} \,\, (\textrm{with  the  index  notation})}
\eeq
Using the definition (\ref{eq: mode-n prod def}) allows to write (\ref{eq: Tucker model a}) in terms of mode-$n$ products as
\beq \label{eq: Tucker model c}
\matc X &=& \matc G {\times}_1 \mat A^{(1)}{\times}_2 \mat A^{(2)}{\times}_3\cdots{\times}_N\mat A^{(N)}\nonumber \\
&=&\matc G{{\times}^{N}_{n=1}}\mat A^{(n)}.
\eeq
This expression evidences that the Tucker model can be viewed as the transformation of the core tensor resulting from its multiplication by the factor matrix $\mat A^{(n)}$ along its mode-$n$, which corresponds to a linear map applied to the mode-$n$ space of $\matc G$, for $n=1,\cdots,N$, i.e. a multilinear map applied to $\matc G$. From a transformation point of view, $\matc G$ and $\matc X$ can be interpreted as the input tensor and the transformed tensor, or output tensor, respectively.

\subsubsection*{\bf{Matrix representations of the Tucker model}}

A matrix representation of a Tucker model is directly linked with a matricization of tensor like (\ref {eq: matricization formula}), corresponding to the combination of two sets of modes $\mathds{S}_1$ and $\mathds{S}_2$. These combinations must be applied both to the tensor $\matc X$ and its core tensor $\matc G$.

The matrix representation (\ref{eq: matricization formula}) of the Tucker model (\ref{eq: Tucker model a}) is given by
\beq \label{eq:Tucker mat rep}
&&{\mat X_{\mathds{S}_1;\mathds{S}_2}=\left( \underset{n \in \mathds{S}_1}{\otimes} \mat A^{(n)} \right)\mat G_{\mathds{S}_1;\mathds{S}_2}\left( \underset{n \in \mathds{S}_2}{\otimes} \mat A^{(n)} \right)^T}
\eeq
with $\mat G_{\mathds{S}_1;\mathds{S}_2}\inmat{J_1}{J_2}$, and $J_{n_1}=\underset{n \in \mathds{S}_{n_1}}{\prod R_n}$, for $n_1=1 \, \textrm{and}\,\,  2$.

\begin{proof}
See the Appendix.
\end{proof}

For the flat mode-$n$ unfolding, defined in (\ref {eq:mode-n unfolding}), the formula (\ref{eq:Tucker mat rep}) gives
\beq \label{eq:mode-n unfolding Tucker}
\mat X_n=\mat A^{(n)}\mat G_n (\mat A^{(n+1)}\otimes \cdots \otimes \mat A^{(N)}\otimes\mat A^{(1)} \otimes \cdots \otimes \mat A^{(n-1)})^T.
\eeq
Applying the vec formula (\ref{eq:vec formula a}) to the right hand-side of (\ref{eq:mode-n unfolding Tucker}), we obtain the vectorized form of $\matc X$ associated with its mode-$n$ unfolding $\mat X_n$
\beq \label{eq:Tucker vectorized form}
\textrm{vec}(\matc X)=\textrm{vec}(\mat X_n)=(\mat A^{(n+1)}\otimes \cdots \otimes \mat A^{(N)}\otimes\mat A^{(1)} \otimes \cdots \otimes \mat A^{(n)})\textrm{vec}(\mat G_n).\nonumber
\eeq

\subsubsection{\bf{Tucker-($N_1,N$) models}}

A Tucker-$(N_1,N)$ model for a $N^{th}$-order tensor $\matc X \inten{I_{1}}{\cdots}{I_{N}}$, with $N\geq N_1$, corresponds to the case where $N-N_1$ factor matrices are equal to identity matrices. For instance, assuming that $\mat A^{(n)}=\mat I_{I_n}$, which implies $R_n=I_n$, for $n=N_1+1,\cdots,N$, Eq. (\ref{eq: Tucker model a}) and (\ref{eq: Tucker model c}) become
\beq \label{eq:Tucker-N1}
x_{i_1,\cdots,i_N }&=&\sum\limits_{r_1=1}^{R_1}\cdots \sum\limits_{r_{N_1}=1}^{R_{N_1}}g_{r_1,\cdots,r_{N_1},i_{N_1+1},\cdots,i_N}\prod\limits_{n=1}^{N_1}a^{(n)}_{i_n,r_n}\\
\label{eq:Tucker-N1 a}
\matc X &=& \matc G {\times}_1 \mat A^{(1)}{\times}_2 \cdots{\times}_{N_1}\mat A^{(N_1)}{\times}_{N_1+1}\mat I_{I_{N_1+1}}\cdots{\times}_N\mat I_{I_N}\nonumber\\
&=&{\matc G}\,{\times}^{N_1}_{n=1} \mat A^{(n)}. \label{eq:Tucker-N1 b}
\eeq
One such model that is currently used in applications is the Tucker-(2,3) model, usually denoted Tucker2, for third-order tensors $\matc X\inten{I}{J}{K}$. Assuming $\mat A^{(1)}=\mat A\inmat{I}{P},\mat A^{(2)}=\mat B\inmat{J}{Q}$, and $\mat A^{(3)}=\mat I_K$, such a model is defined by the following equations
\beq \label{eq:Tucker2}
x_{ijk}&=&\sum\limits_{p=1}^{P}\sum\limits_{q=1}^{Q}g_{pqk}a_{ip}b_{jq}\\
\matc X&=&\matc G{\times}_1 \mat A{\times}_2 \mat B\label{eq:Tucker2 a}
\eeq
with the core tensor $\matc G\inten{P}{Q}{K}$.

\subsubsection{\bf{PARAFAC Models}}
A PARAFAC model for a $N^{th}$-order tensor corresponds to the particular case of a Tucker model with an identity core tensor of order $N$ and dimensions $R\times\cdots\times R$
\beq
&&{\matc G={\matc I}_{N,R}=\matc I} \quad \Leftrightarrow \quad g_{r_1,\cdots,r_N }=\delta_{r_1,\cdots,r_N }\nonumber
\eeq
Equations (\ref{eq: Tucker model a})-(\ref{eq: Tucker model c}) then become, respectively
\beq \label {eq: Parafac model a}
x_{i_1,\cdots,i_N }&=&\sum\limits_{r=1}^{R}\prod\limits_{n=1}^{N}a^{(n)}_{i_n,r}\\
&=&{\color{black}\prod\limits_{n=1}^{N}a^{(n)}_{i_n,r}\,\,\,\,\,\,(\textrm{with the index notation})}\label{PARAFAC ind not}\\
\label {eq: Parafac model b}
\matc X &=& \sum\limits_{r=1}^{R}{(\overset{N}{\underset{n=1}{\circ}}\mat A^{(n)}_{.r})}\nonumber\\
\label {eq: Parafac model c}
\matc X&=&\matc I_{N,R} {\times}^{N}_{n=1} \mat A^{(n)}
\eeq
with the factor matrices $\mat A^{(n)}\inmat{I_n}{R},n=1,\cdots,N$.

\subsection*{Remarks}

\begin{itemize}
\item The expression (\ref{eq: Parafac model a}) as a sum of polyads is called a polyadic form of $\matc X$ by Hitchcock (1927) \cite{Hitchcock27}.
\item The PARAFAC model (\ref{eq: Parafac model a})-(\ref{eq: Parafac model c}) amounts to decomposing the tensor $\matc X$ into a sum of $R$ components, each component being a rank-one tensor. When $R$ is minimal in (\ref {eq: Parafac model a}), it is called the rank of $\matc X$ \cite{Kru77}. This rank is related to the mode-$n$ ranks by the following inequalities $\textrm{rank}_n(\matc X)\leq R, \forall n=1, \cdots, N$. Furthermore, contrary to the matrices for which the rank is always at most equal to the smallest of the dimensions, for higher-order tensors the rank can exceed any mode-$n$ dimension $I_n$.

    There exists different definitions of rank for tensors, like typical and generic ranks, or also symmetric rank for a symmetric tensor. See \cite{Com09:LinearAlgebra} and \cite{ComoGLM08:SIAM} for more details.
\item In telecommunication applications, the structure parameters (rank, mode dimensions, and core tensor dimensions) of a PARAFAC or Tucker model, are design parameters that are chosen in function of the performance desired for the communication system. However, in most of the applications, as for instance in multi-way data analysis, the structure parameters are generally unknown and must be determined a priori. Several techniques have been proposed for determining these parameters. See \cite{Bro2003}, \cite{daCosta2008}, \cite{daCosta2010}, \cite{daCosta2011}, and references therein.
\item The PARAFAC model is also sometimes defined by the following equation
\beq \label {eq: Parafac model d}
&&{x_{i_1,\cdots,i_N }=\sum\limits_{r=1}^{R}g_r\prod\limits_{n=1}^{N}a^{(n)}_{i_n,r} \quad \textrm{with} \quad g_r>0}.
\eeq

In this case, the identity tensor $\matc I_{N,R}$ in (\ref{eq: Parafac model c}) is replaced by the diagonal tensor $\matc G \inten {R}{\cdots}{R}$ whose diagonal elements are equal to scaling factors $g_r$, i.e.
\
\beq
g_{r_1,\cdots,r_N }=\left\{\begin{array}{ll} g_r \quad \textrm{if}\quad r_1=\cdots=r_N=r\\
0 \quad \textrm{otherwise} \end{array}\right.\nonumber
\eeq
and all the column vectors $\mat A^{(n)}_{.r}$ are normalized, i.e. with a unit norm, for $1\leq n\leq N$.
\end{itemize}
\begin{itemize}
\item It is important to notice that the PARAFAC model (\ref {eq: Parafac model a}) is multilinear (more precisely $N$-linear) in its parameters in the sense that it is linear with respect to each matrix factor. This multilinearity property is exploited for parameter estimation using the standard alternating least squares (ALS) algorithm (\cite {Harshman70}, \cite {CarrollChang70}) that consists in alternately estimating each matrix factor by minimizing a least squares error criterion conditionally to the knowledge of the other matrix factors that are fixed with their previously estimated values.

\end{itemize}

\subsubsection*{\bf{Matrix representations of the PARAFAC model}}

The matrix representation (\ref{eq: matricization formula}) of the PARAFAC model (\ref{eq: Parafac model a})-(\ref {eq: Parafac model c}) is given by
\beq \label{eq:Parafac mat rep}
\mat X_{\mathds{S}_1;\mathds{S}_2}=\left( \underset{n \in \mathds{S}_1}{\khatri} \mat A^{(n)} \right)\left( \underset{n \in \mathds{S}_2}{\khatri} \mat A^{(n)} \right)^T.
\eeq

\begin{proof}
See the Appendix.
\end{proof}

\subsection*{Remarks}

\begin{itemize}
\item From (\ref {eq:Parafac mat rep}), we can deduce that a mode combination results in a Khatri-Rao product of the corresponding factor matrices. Consequently, the tensor contraction (\ref{eq: contraction def}) associated with the PARAFAC-$N$ model (\ref {eq: Parafac model c}) gives a PARAFAC-$N_1$ model whose factor matrices are equal to $\underset{n \in \mathds{S}_{n_1}}{\khatri} \mat A^{(n)}\inmat{J_{n_1}}{R}$, $n_1=1,\cdots,N_1$, with $J_{n_1}=\underset{n \in \mathds{S}_{n_1}}{\prod I_n}$.

\item For the PARAFAC model, the flat mode-$n$ unfolding, defined in (\ref {eq:mode-n unfolding}), is given by
\beq \label{eq:mode-n unfolding PARAFAC}
\mat X_n=\mat A^{(n)}(\mat A^{(n+1)}\khatri \cdots \khatri \mat A^{(N)}\khatri\mat A^{(1)} \khatri \cdots \khatri \mat A^{(n-1)})^T
\eeq
and the associated vectorized form is obtained in applying the vec formula (\ref{eq:vec formula b}) to the right hand-side of the above equation, with $\mat I_R=\textrm{diag}(\mat 1_R)$
\beq \label{eq:PARAFAC vectorized form}
\textrm{vec}(\matc X)=\textrm{vec}(\mat X_n)=(\mat A^{(n+1)}\khatri \cdots \khatri \mat A^{(N)}\khatri\mat A^{(1)} \khatri \cdots \khatri \mat A^{(n)})\mat 1_{R}
\eeq

\item In the case of the normalized PARAFAC model (\ref {eq: Parafac model d}), Eq. (\ref {eq:Parafac mat rep}) and (\ref{eq:PARAFAC vectorized form}) become, respectively
\beq \label{eq:normalized PARAFAC mat rep}
\mat X_{\mathds{S}_1;\mathds{S}_2}&=&\left( \underset{n \in \mathds{S}_1}{\khatri} \mat A^{(n)} \right)\textrm{diag}(\mat g)\left( \underset{n \in \mathds{S}_2}{\khatri} \mat A^{(n)} \right)^T\nonumber\\
\label{eq:normalized PARAFAC vectorized form}
\textrm{vec}(\matc X)&=&\textrm{vec}(\mat X_n)=(\mat A^{(n+1)}\khatri \cdots \khatri \mat A^{(N)}\khatri\mat A^{(1)} \khatri \cdots \khatri \mat A^{(n)})\mat g\nonumber
\eeq
where $\mat g=[g_1 \cdots g_R]^T\inmat{R}{1}$.



\item For the PARAFAC model of a third-order tensor $\matc X \inten{I}{J}{K}$ with factor matrices $(\mat A,\mat B,\mat C)$, the formula (\ref{eq:Parafac mat rep}) gives for $\mathds{S}_1=\{i,j\}$ and $\mathds{S}_2=\{k\}$
\beq
\label{eq:third order unfoldings a}
\mat X_{IJ\times K}=\left[\begin{array}{c}\mat X_{1..}\\ \vdots \\ \mat X_{I..}\end{array}\right]=(\mat A \khatri \mat B)\mat C^T \inmat{IJ}{K}.\nonumber
\eeq
Noting that $\mat A \khatri \mat B=\left[\begin{array}{c}\mat BD_1(\mat A)\\ \vdots \\ \mat BD_I(\mat A)\end{array}\right]$, we deduce the following expression for mode-1 matrix slices
\vspace{-1ex}
\beq \label{eq: third order slices}
\mat X_{i..}=\mat B D_i(\mat A)\mat C^T.\nonumber
\eeq
Similarly, we have
\beq
\mat X_{JK\times I}&=&(\mat B \khatri \mat C)\mat A^T, \quad \mat X_{KI\times J}=(\mat C \khatri \mat A)\mat B^T,\nonumber\\
\mat X_{.j.}&=&\mat C D_j(\mat B)\mat A^T,\quad \mat X_{..k}=\mat A D_k(\mat C)\mat B^T.\nonumber
\eeq

\end{itemize}

\begin{itemize}
\item For the PARAFAC model of a fourth-order tensor $\matc X \inten{I}{J}{K\times L}$ with factor matrices $(\mat A,\mat B,\mat C, \mat D)$, we obtain
\beq
\mat X_{IJK\times L}&=&(\mat A\khatri\mat B\khatri\mat C)\mat D^T\nonumber\\
&=&\left[\begin{array}{c}(\mat B\khatri\mat C)D_1(\mat A)\\ \vdots \\ (\mat B\khatri\mat C)D_I(\mat A)\end{array}\right]\mat D^T = \left[\begin{array}{c}\mat CD_1(\mat B)D_1(\mat A)\\ \vdots \\ \mat CD_J(\mat B)D_I(\mat A)\end{array}\right]\mat D^T \inmat{IJK}{L} \nonumber\\
\label{eq: fourth order slices}
\mat X_{ij..}&=&\mat C D_j(\mat B)D_i(\mat A)\mat D^T\inmat{K}{L}
\eeq
Other matrix slices can be deduced from (\ref{eq: fourth order slices}) by simple permutations of the matrix factors.
\end{itemize}

In the next section, we introduce two constrained PARAFAC models, the so called PARALIND and CONFAC models, and then PARATUCK models.

{\color{black}
\section{\bf{Constrained PARAFAC Models}}}

The introduction of constraints in tensor models can result from the system itself that is under study, or from a system design.
In the first case, the constraints are often interpreted as interactions or linear dependencies between the PARAFAC factors. Examples of such dependencies are encountered in psychometrics and chemometrics applications that gave origin, respectively, to the PARATUCK-2 model \cite{Harshman96} and the PARALIND (PARAllel profiles with LINear Dependencies) model (\cite{Bro05}, \cite{BHSL09}), introduced in \cite{CPK80} under the name CANDELINC (CANonical DEcomposition with LINear Constraints), for the multiway case. A first application of the PARATUCK-2 model in signal processing was made in \cite{Kibangou07Eusipco} for blind joint identification and equalization of Wiener-Hammerstein communication channels. The PARALIND model was recently applied for identifiability and propagation parameter estimation purposes in a context of array signal processing \cite{Xu2012a}, \cite{Xu2012b}.

In the second case, the constraints are used as design parameters. For instance, in a telecommunications context, we recently proposed two constrained tensor models: the CONFAC (CONstrained FACtor) model \cite{AFM08IEEESP}, and the PARATUCK-$(N_1,N)$ model \cite{FCAR11EUSIPCO}, \cite{FCAR12SP}. The PARATUCK-2 model was also applied for designing space-time spreading-multiplexing MIMO systems \cite{Almeida2009}. For these telecommunication applications of constrained tensor models, the constraints are used for resource allocation. We are now going to describe these various constrained PARAFAC models.\\

\subsection{\bf{PARALIND models}}

Let us define the core tensor of the Tucker model (\ref {eq: Tucker model c}) as follows:
\beq \label{eq: Paralind core}
\matc G={\matc I}_{N,R}\,{\times}^{N}_{n=1} \mat \Phi^{(n)}
\eeq
where  $\mat \Phi^{(n)}\inmatR{R_n}{R}, \, n = 1,\cdots, N$, with $R\geq \underset{n}{max}(R_n)$, are constraint matrices. In this case, $\matc G$ will be called the "interaction tensor", or "constraint tensor".

The PARALIND model is obtained by substituting (\ref{eq: Paralind core}) into (\ref{eq: Tucker model c}), and applying the property (\ref{eq: mode-n prod property}), which gives
\vspace{-3ex}
\beq \label {eq: Paralind model a}
\matc X=\matc G {\times}^{N}_{n=1} \mat A^{(n)}={\matc I}_{N,R}\,{\times}^{N}_{n=1} (\mat A^{(n)}\mat \Phi^{(n)}).
\eeq
Equation (\ref{eq: Paralind model a}) leads to two different interpretations of the PARALIND model, as a constrained Tucker model whose core tensor admits a PARAFAC decomposition with factor matrices  $\mat \Phi^{(n)}$, called "interaction matrices", and as a constrained PARAFAC model with constrained factor matrices $\bar{\mat A}^{(n)}=\mat A^{(n)}\mat \Phi^{(n)}$.

The interaction matrix $\mat \Phi^{(n)}$ allows taking into account linear dependencies between the columns of  $\mat A^{(n)}$, implying a rank deficiency for this factor matrix. When the columns of  $\mat \Phi^{(n)}$ are formed with $0's$ and $1's$, the dependencies simply consist in a repetition or an addition of certain columns of $\mat A^{(n)}$. In this particular case, the diagonal element $\xi^{(n)}_{r,r}\geq1$ of the matrix $\mat \Xi^{(n)}={\mat \Phi^{(n)}}^T \mat \Phi^{(n)}\inmatR{R}{R}$, represents the number of columns of $\mat A^{(n)}$  that are added to form the $r^{th}$ column of the constrained factor $\mat A^{(n)}\mat \Phi^{(n)}$. The choice $\mat \Phi^{(n)}=\mat I_R$  means that there is no such dependency among the columns of $\mat A^{(n)}$.

Equation (\ref{eq: Paralind model a}) can be written element-wise as
\beq
x_{i_1,\cdots,i_N }&=&\sum\limits_{r_1=1}^{R_1}\cdots \sum\limits_{r_N=1}^{R_N}g_{r_1,\cdots,r_N }\prod\limits_{n=1}^{N}a^{(n)}_{i_n,r_n} \,\, \textrm{with} \,\, g_{r_1,\cdots,r_N }=\sum\limits_{r=1}^{R}\prod\limits_{n=1}^{N}\phi^{(n)}_{r_n,r}\nonumber\\
\label{eq: Paralind model c}
&=&\sum\limits_{r=1}^{R}\prod\limits_{n=1}^{N}\bar{a}^{(n)}_{i_n,r} \,\, \textrm{with} \,\, \bar{a}^{(n)}_{i_n,r}=\sum\limits_{r_n=1}^{R_n}a^{(n)}_{i_n,r_n}\phi^{(n)}_{r_n,r}.
\eeq
This constrained PARAFAC model constitutes an $N$-way form of the three-way PARALIND model, used for chemometrics applications in \cite{Bro05}, and \cite{BHSL09}.

\subsection{\bf{CONFAC models}}

When the constraint matrices $\mat \Phi^{(n)}\inmatR{R_n}{R}$ are full row-rank, and their columns are chosen as canonical vectors of the Euclidean space $\mathbb{R}^{R_n}$, for $n = 1,\cdots, N$, the constrained PARAFAC model (\ref{eq: Paralind model a}) constitutes a generalization to $N^{th}$-order of the third-order CONFAC model, introduced in \cite{AFM08IEEESP} for designing MIMO communication systems with resource allocation. {\color{black} This CONFAC model was used in \cite{ALSC2012} for solving the problem of blind identification of underdetermined mixtures based on cumulant generating function of the observations.}
In a telecommunications context where $\matc X$ represents the tensor of received signals, such a constraint matrix $\mat \Phi^{(n)}$ can be interpreted as an "allocation matrix" allowing to allocate resources, like data streams, codes, and transmit antennas, to the $R$ components of the signal to be transmitted. In this case, the core tensor $\matc G$ will be called the "allocation tensor". By assumption, each column of the allocation matrix $\mat \Phi^{(n)}$ is a canonical vector of $\mathbb{R}^{R_n}$, which means that there is only one value of $r_n$ such that $\phi^{(n)}_{r_n,r}=1$, and this value of $r_n$  corresponds to the $n^{th}$ resource allocated to the $r^{th}$ component.

Each element $x_{i_1,\cdots,i_N }$ of the received signal tensor $\matc X$ is equal to the sum of $R$ components, each component $r$ resulting from the combination of $N$ resources, each resource being associated with a column of the matrix factor $\mat A^{(n)}$, $n = 1,\cdots, N$. This combination, determined by the allocation matrices, is defined by a set of $N$ indices $\{r_1,\cdots,r_N\}$ such that $\prod\limits_{n=1}^{N}\phi^{(n)}_{r_n,r}=1$. As for any $r\in[1,R]$, there is one and only one $N$-uplet $(r_1,\cdots ,r_N)$ such as $\prod\limits_{n=1}^N\phi_{r_n,r}^{(n)}=1$, we can deduce that each component $r$ of $x_{i_1,\cdots,i_N }$  in (\ref{eq: Paralind model c}) is the result of one and only one combination of the $N$ resources under the form of the product $\prod\limits_{n=1}^{N}a^{(n)}_{i_n,r_n}$.
For the CONFAC model, we have
\beq
\sum\limits_{r_n=1}^{R_n}D_{r_n}(\mat \Phi^{(n)})=\mat I_R,\,\, \forall n=1,\cdots ,N\nonumber
\eeq
meaning that each resource $r_n$ is allocated at least once, and the diagonal element of $\mat \Xi^{(n)}={\mat \Phi^{(n)}}^T \mat \Phi^{(n)}$ is such as $\xi^{(n)}_{r,r}=1, \forall n=1,\cdots ,N$, because only one resource $r_n$ is allocated to each component $r$. Moreover, we have to notice that the assumption $R\geq\underset{n}{max}(R_n)$ implies that each resource can be allocated several times, i.e. to several components. Defining the interaction matrices
\beq
\mat \Gamma^{(n)}=\mat \Phi^{(n)} {\mat \Phi^{(n)}}^T \inmatR{R_n}{R_n}, \mat \Gamma^{(n_1,n_2)}=\mat \Phi^{(n_1)} {\mat \Phi^{(n_2)}}^T \inmatR{R_{n_1}}{R_{n_2}}\nonumber
\eeq
the diagonal element $\gamma^{(n)}_{r_n,r_n}\in[1,R-R_n+1]$ represents the number of times that the $r^{th}_{n}$ column of $\mat A^{(n)}$ is repeated, i.e. the number of times that the $r^{th}_{n}$  resource is allocated to the $R$ components, whereas  $\gamma^{(n_1,n_2)}_{r_{n_1},r_{n_2}}$ determines the number of interactions between the $r^{th}_{n_1}$ column of $\mat A^{(n_1)}$ and the $r^{th}_{n_2}$ column of $\mat A^{(n_2)}$, i.e. the number of times that the $r^{th}_{n_1}$ and $r^{th}_{n_2}$ resources are combined in the $R$ components. If we choose $R_n=R$ and $\mat \Phi^{(n)}=\mat I_R, \forall n=1,\cdots,N$, the PARALIND/CONFAC model (\ref{eq: Paralind model a}) becomes identical to the PARAFAC one (\ref{eq: Parafac model c}).

The matrix representation (\ref{eq: matricization formula}) of the PARALIND/CONFAC model can be deduced from (\ref{eq:Parafac mat rep}) in replacing $\mat A^{(n)}$ by $\mat A^{(n)}\mat \Phi^{(n)}$
\beq
\mat X_{\mathds{S}_1;\mathds{S}_2}=\left( \underset{n \in \mathds{S}_1}{\khatri} \mat A^{(n)}\mat \Phi^{(n)} \right)\left( \underset{n \in \mathds{S}_2}{\khatri} \mat A^{(n)}\mat \Phi^{(n)} \right)^T.\nonumber
\eeq
Using the identity (\ref {eq: kron prop3}) gives
\beq \label{eq:CONFAC mat rep a}
\mat X_{\mathds{S}_1;\mathds{S}_2}=\left( \underset{n \in \mathds{S}_1}{\kron} \mat A^{(n)} \right)\left( \underset{n \in \mathds{S}_1}{\khatri} \mat \Phi^{(n)} \right)\left( \underset{n \in \mathds{S}_2}{\khatri} \mat \Phi^{(n)} \right)^T\left( \underset{n \in \mathds{S}_2}{\kron} \mat A^{(n)} \right)^T,
\eeq
or, equivalently,
\beq \label{eq:CONFAC mat rep b}
\mat X_{\mathds{S}_1;\mathds{S}_2}=\left( \underset{n \in \mathds{S}_1}{\kron} \mat A^{(n)} \right)\mat G_{\mathds{S}_1;\mathds{S}_2}\left( \underset{n \in \mathds{S}_2}{\kron} \mat A^{(n)} \right)^T,\nonumber
\eeq
where the matrix representation $\mat G_{\mathds{S}_1;\mathds{S}_2}$ of the constraint/allocation tensor $\matc G$, defined by means of its PARAFAC model (\ref {eq: Paralind core}), can also be deduced from (\ref{eq:Parafac mat rep}) as
\beq \label{eq:CONFAC mat rep c}
\mat G_{\mathds{S}_1;\mathds{S}_2}=\left( \underset{n \in \mathds{S}_1}{\khatri} \mat \Phi^{(n)} \right)\left( \underset{n \in \mathds{S}_2}{\khatri} \mat \Phi^{(n)} \right)^T.\nonumber
\eeq

\subsection{\bf{Nested Tucker models}}

The PARALIND/CONFAC models can be viewed as particular cases of a new family of tensor models that we shall call nested Tucker models, defined by means of the following recursive equation
\beq \label {eq: nested Tucker models}
\matc X^{(p)} &=& \matc X^{(p-1)}{\times}^{N}_{n=1}\mat A^{(p,n)} \,\,\, \textrm{for} \,\,\, p=1,\cdots, P \nonumber\\
&=& \matc G {\times}^{N}_{n=1}\prod\limits_{q=P}^{1}\mat A^{(q,n)}\nonumber
\eeq
\
with the factor matrices $\mat A^{(p,n)}\inmat{R^{(p,n)}}{R^{(p-1,n)}}$ for $p=1,\cdots, P$, such as $R^{(0,n)}=R_n$ and $R^{(P,n)}=I_n$, for $n=1,\cdots,N$, the core tensor $\matc X^{(0)}=\matc G \inten{R_{1}}{\cdots}{R_{N}}$, and $\matc X^{(P)}\inten{I_{1}}{\cdots}{I_{N}}$. This equation can be interpreted as $P$ successive linear transformations applied to each mode-$n$ space of the core tensor $\matc G$. So, $P$ nested Tucker models can then be interpreted as a Tucker model for which the factor matrices are products of $P$ matrices. When $\matc G = \matc I_{N,R}$, which implies $R^{(0,n)}=R_n=R$ for $n=1,\cdots,N$, we obtain nested PARAFAC models. The PARALIND/CONFAC models correspond to two nested PARAFAC models ($P=2$), with $\mat A^{(1,n)}=\mat \Phi^{(n)}$, $\mat A^{(2,n)}=\mat A^{(n)}$, $R^{(0,n)}=R$, $R^{(1,n)}=R_n$, and $R^{(2,n)}=I_n$, for $n=1,\cdots,N$.

By considering nested PARAFAC models with $P=3$, $\mat A^{(1,n)}=\mat \Phi^{(n)}\inmat{K_n}{R}$, $\mat A^{(2,n)}=\mat A^{(n)}\inmat{J_n}{K_n}$ and $\mat A^{(3,n)}=\mat \Psi^{(n)}\inmat{I_n}{J_n}$, for $n=1,\cdots,N$, we deduce doubly PARALIND/CONFAC models described by the following equation
\beq \label {eq: doubly Paralind model}
\matc X={\matc I}_{N,R}\,{\times}^{N}_{n=1} (\mat \Psi^{(n)}\mat A^{(n)}\mat \Phi^{(n)}).\nonumber
\eeq
Such a model can be viewed as a doubly constrained PARAFAC model, with factor matrices $\mat \Psi^{(n)}\mat A^{(n)}\mat \Phi^{(n)}$, the constraint matrix $\mat \Psi^{(n)}$, assumed to be full column-rank, allowing to take into account linear dependencies between the rows of $\mat A^{(n)}$.

\begin{figure}[!t]
\centering
\includegraphics[width=.99\columnwidth]{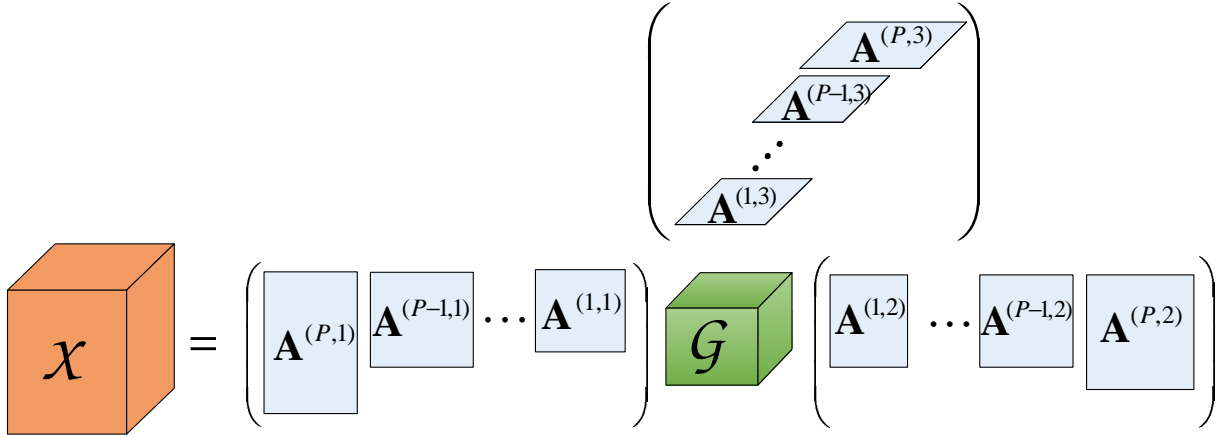}
\caption{Visualization of the nested Tucker model.} \label{fig:1}
\end{figure}

\subsection{\bf{Block PARALIND/CONFAC models}}

In some applications, the data tensor $\matc X \inten{I_{1}}{\cdots}{I_{N}}$ is written as a sum of $P$ sub-tensors $\matc X^{(p)}$, each sub-tensor admitting a tensor model with a possibly different structure. So, we can define a block-PARALIND/CONFAC model as
\beq \label {eq: block tensor model}
\matc X&=&\sum\limits_{p=1}^{P}\matc X^{(p)},\\
\label {eq: block PARALIND/CONFAC model}
\matc X^{(p)}&=&\matc G^{(p)} {\times}^{N}_{n=1} \mat A^{(p,n)},\\
\matc G^{(p)}&=&{\matc I}_{N,R^{(p)}}\,{\times}^{N}_{n=1} \mat \Phi^{(p,n)},\nonumber
\eeq
where $\mat A^{(p,n)}\inmat{I_n}{R^{(p,n)}}$, $\mat \Phi^{(p,n)}\inmat{R^{(p,n)}}{R^{(p)}}$, and $\matc G^{(p)}\inten{R^{(p,1)}}{\cdots}{R^{(p,N)}}$ are the mode-$n$ factor matrix, the mode-$n$ constraint/allocation matrix, and the core tensor of the PARALIND/CONFAC model of the $p^{th}$ sub-tensor, respectively. The matrix representation (\ref{eq:CONFAC mat rep a}) then becomes
\beq \label{eq:Block CONFAC mat rep a}
\mat X_{\mathds{S}_1;\mathds{S}_2}=\sum\limits_{p=1}^{P}\left( \underset{n \in \mathds{S}_1}{\kron} \mat A^{(p,n)} \right)\left( \underset{n \in \mathds{S}_1}{\khatri} \mat \Phi^{(p,n)} \right)\left( \underset{n \in \mathds{S}_2}{\khatri} \mat \Phi^{(p,n)} \right)^T\left( \underset{n \in \mathds{S}_2}{\kron} \mat A^{(p,n)} \right)^T.
\eeq
Defining the following block partitioned matrices
\beq \label{eq:matrix factor partitioning}
\mat A^{(n)}=\left[\begin{array}{ccc}\mat A^{(1,n)} \cdots  \mat A^{(P,n)}\end{array}\right]\inmat {I_n}{R^{(n)}}
\eeq
where $R^{(n)}=\sum\limits_{p=1}^{P}R^{(p,n)}$, Eq. (\ref{eq:Block CONFAC mat rep a}) can be rewritten in the following more compact form
\vspace{-1ex}
\beq \label{eq:Block CONFAC mat rep b}
\mat X_{\mathds{S}_1;\mathds{S}_2}=\left( \underset{n \in \mathds{S}_1}{{\otimes}_b} \mat A^{(n)} \right)\mat G_{\mathds{S}_1;\mathds{S}_2}\left( \underset{n \in \mathds{S}_2}{{\otimes}_b} \mat A^{(n)} \right)^T\nonumber
\eeq
where ${\otimes}_b$ denotes the block-wise Kronecker product defined as
\beq \label {eq:Block Kro}
\mat A^{(n)}{\otimes}_b \mat A^{(q)}=\left[\begin{array}{ccc}\mat A^{(1,n)}\otimes\mat A^{(1,q)} \cdots  \mat A^{(P,n)}\otimes\mat A^{(P,q)}\end{array}\right]\nonumber
\eeq
$\mat A^{(q)}$ being partitioned in $P$ blocks as in (\ref{eq:matrix factor partitioning}), and
\beq
\mat G_{\mathds{S}_1;\mathds{S}_2}&=&bdiag(\mat G_{\mathds{S}_1;\mathds{S}_2}^{(1)}\cdots \mat G_{\mathds{S}_1;\mathds{S}_2}^{(P)})\inmat{J_1}{J_2}\nonumber\\
\mat G_{\mathds{S}_1;\mathds{S}_2}^{(p)}&=&\left( \underset{n \in \mathds{S}_1}{{\khatri}_b} \mat \Phi^{(p,n)} \right)\left( \underset{n \in \mathds{S}_2}{{\khatri}_b} \mat \Phi^{(p,n)} \right)^T\inmat{J_1^{(p)}}{J_2^{(p)}}\nonumber
\eeq
where ${\khatri}_b$ denotes the block-wise Khatri-Rao product defined in the same way as the block-wise Kronecker product, with $J_{n_1}=\sum\limits_{p=1}^{P}{J_{n_1}^{(p)}}$ and $J_{n_1}^{(p)}=\underset{n \in \mathds{S}_{n_1}}{\prod}R^{(p,n)}$ for $n_1=1 \, \textrm{and}\,\,  2$.
\\

In the case of a block PARAFAC model, Eq. (\ref {eq: block PARALIND/CONFAC model}) is replaced by
\beq \label {eq: block PARAFAC model}
\matc X^{(p)}=\matc I_{N,R^{(p)}} {\times}^{N}_{n=1} \mat A^{(p,n)} \,\,\,  \textrm{with} \,\,\,  \mat A^{(p,n)}\inmat{I_n}{R^{(p)}}\nonumber
\eeq
and the matrix representation (\ref {eq:Parafac mat rep}) then becomes
\beq \label {eq:Block Parafac mat rep}
\mat X_{\mathds{S}_1;\mathds{S}_2}=\left( \underset{n \in \mathds{S}_1}{{\khatri}_b} \mat A^{(n)} \right)\left( \underset{n \in \mathds{S}_2}{{\khatri}_b} \mat A^{(n)} \right)^T\nonumber
\eeq
\
with $\mat A^{(n)}=\left[\begin{array}{ccc}\mat A^{(1,n)} \cdots  \mat A^{(P,n)}\end{array}\right]\inmat {I_n}{R}$, and $R=\sum\limits_{p=1}^{P}{R^{(p)}}$. Block constrained PARAFAC models were used in \cite{Almeida05Asilomar}, \cite{Almeida05PSIP}, \cite{Almeida07SP} for modeling different types of multiuser wireless communication systems. Block constrained Tucker models were used for space-time multiplexing MIMO-OFDM systems \cite{Almeida06PIMRC}, and for blind beamforming \cite{Almeida2009b}. In these applications, the symbol matrix factor is in Toeplitz or block-Toeplitz form.

The block tensor model defined by Eq. (\ref {eq: block tensor model})-(\ref {eq: block PARALIND/CONFAC model}) can be viewed as a generalization of the block term decomposition introduced in \cite{DeLathauwer2008} for third-order tensors $\matc X\inten{I}{J}{K}$ that are decomposed into a sum of $P$ Tucker models of rank-$(L,M,N)$, which corresponds to the particular case where all the factor matrices are full column rank, with $\mat A^{(p,1)}\inmat{I}{L}$, $\mat A^{(p,2)}\inmat{J}{M}$, and $\mat A^{(p,3)}\inmat{K}{N}$, for $p=1,\cdots,P$, and $\matc G\inten{L}{M}{N}$, and each sub-tensor $\matc X^{(p)}$ is decomposed by means of its HOSVD.

\begin{figure}[!t]
\centering
\includegraphics[width=.99\columnwidth]{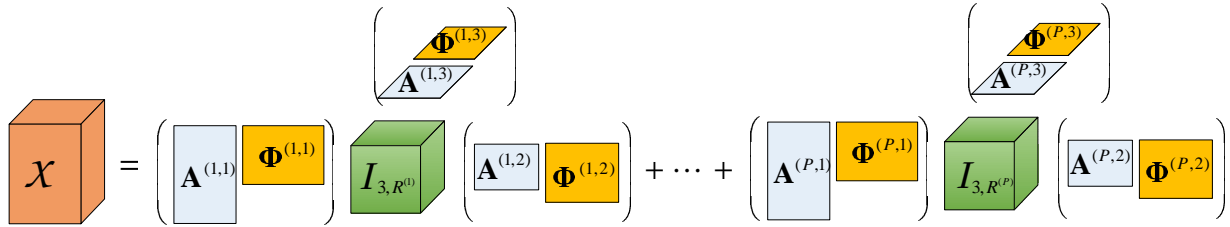}
\caption{Visualization of the block PARALIND/CONFAC model.} \label{fig:2}
\end{figure}

{\color{black} This figure is to be compared with Figure 5 in \cite{CMP2014} representing a block term decomposition of a third-order tensor into rank-$(L_p,M_p,N_p)$ terms, when each term has a PARALIND/CONFAC structure.}

\subsection{\bf{PARALIND/CONFAC-$(N_1,N)$ models}}

Now, we introduce a variant of PARALIND/CONFAC models that we shall call PARALIND/CONFAC-$(N_1,N)$ models. This variant corresponds to PARALIND/CONFAC models (\ref{eq: Paralind model a}) with only $N_1$ constrained matrix factors, which implies $R_n=R$ and $\mat A^{(n)}\inmat{I_n}{R}$ for $n=N_1+1,\cdots ,N$
\beq \label {eq: Paralind N_1 model a}
\matc X={\matc I}_{N,R}\,{\times}^{N_1}_{n=1} (\mat A^{(n)}\mat \Phi^{(n)})\,{\times}^{N}_{n=N_1+1}\mat A^{(n)}.
\eeq
In \cite{AFM08IEEESP:2}, a block PARALIND/CONFAC-(2,3) model that can be deduced from (\ref {eq: Paralind N_1 model a}), was used for modeling uplink multiple-antenna code-division multiple-access (CDMA) multiuser systems.\\

{\color{black}The block term decomposition (BTD) in rank-$(1,L_p,L_p)$ terms of a third-order tensor $\matc X\inten{I}{J}{K}$, which is compared to a third-order PARATREE model in \cite{SRK2009}, can also be viewed as a particular CONFAC-(1,3) model. Indeed, such a decomposition can be written as \cite{L2011}
\beq \label{BTD}
\matc X=\sum\limits_{p=1}^{P}\mat a_p \circ (\mat B_p\mat C_p^T)
\eeq
where the matrices $\mat B_p\inmat{J}{L_p}$ and $\mat C_p\inmat{K}{L_p}$ are rank-$L_p$, and $\mat a_p\inmat{I}{1}$. Defining $\mat B=[\mat B_1\cdots \mat B_P]\inmat{J}{R}$, $\mat C=[\mat C_1\cdots \mat C_P]\inmat{K}{R}$, and $\mat A=[\mat a_1\cdots \mat a_P]\inmat{I}{P}$, with $R=\sum\limits_{p=1}^{P}L_p$, it is easy to verify that the BTD (\ref{BTD}) can be rewritten as the following CONFAC-(1,3) model
\beq
\matc X= \matc I_{3,R}\times_1 \mat A {\bf{\Phi}}\times_2 \mat B\times_3 \mat C
\eeq
with the constraint matrix $\mat \Phi=\left[\begin{array}{ccc} \mat 1_{L_1}^T & & \\
& \ddots &\\
& & \mat 1_{L_P}^T \end{array}\right]\inmat{P}{R}$.
}

\subsection{\bf{PARATUCK models}}

A PARATUCK-$(N_1, N)$ model for a $N^{th}$-order tensor $\matc X \inten{I_{1}}{\cdots}{I_{N}}$, with $N > N_1$, is defined in scalar form as follows \cite{FCAR11EUSIPCO}, \cite{FCAR12SP}
\beq \label{eq:general Paratuck model}
x_{i_1,\cdots,i_{N_1+1},\cdots,i_N }=\sum\limits_{r_1=1}^{R_1}\cdots \sum\limits_{r_{N_1}=1}^{R_{N_1}}c_{r_1,\cdots,r_{N_1},i_{{N_1}+2},\cdots,i_N}\prod\limits_{n=1}^{N_1}a^{(n)}_{i_n,r_n}\phi^{(n)}_{r_n,i_{{N_1}+1}}
\eeq
where $a^{(n)}_{i_n,r_n}$, and $\phi^{(n)}_{r_n,i_{{N_1}+1}}$ are entries of the factor matrix $\mat A^{(n)}\inmat{I_n}{R_n}$ and of the interaction/allocation matrix $\mat \Phi^{(n)}\inmat{R_n}{I_{{N_1}+1}},\, \forall n=1,\cdots,N_1$, respectively, and $\matc C\inten{R_{1}\times\cdots}{R_{N_1}\times I_{N_1+2}\times\cdots}{I_N}$ is the $(N-1)^{th}$-order input tensor. Defining the core tensor $\matc G \inten{R_{1}\times\cdots}{R_{N_1}\times I_{N_1+1}\times\cdots}{I_N}$ element-wise as
\beq \label{eq:general Paratuck core}
g_{r_1,\cdots,r_{N_1},i_{{N_1}+1},\cdots,i_N}&=&c_{r_1,\cdots,r_{N_1},i_{{N_1}+2},\cdots,i_N}\prod\limits_{n=1}^{N_1}\phi^{(n)}_{r_n,i_{{N_1}+1}}\nonumber
\eeq
the PARATUCK-$(N_1,N)$ model can be rewritten as a Tucker-$(N_1,N)$ model (\ref{eq:Tucker-N1})-(\ref{eq:Tucker-N1 b}).


Defining the allocation/interaction tensor $\matc F\inten{R_{1}\times\cdots}{R_{N_1}}{I_{N_1+1}}$ of order $N_1+1$, such as
\beq \label{eq:Paratuck allocation tensor}
f_{r_1,\cdots,r_{N_1},i_{N_1+1}}=\prod\limits_{n=1}^{N_1}\phi^{(n)}_{r_n,i_{N_1+1}}.
\eeq
the core tensor $\matc G$ can then be written as the Hadamard product of the tensors $\matc C$ and $\matc F$ along their first $N_1$ modes
\beq \label{eq:Hadamard products}
\matc G=\,\matc C\, \underset{\{r_1,\cdots,r_{N_1}\}}{\odot}\, \matc F.
\eeq


\subsection*{Remarks}

\begin{itemize}
\item The PARATUCK-$(N1, N)$ model can be interpreted as the transformation of the input tensor $\matc C$  via its multiplication by the factor matrices  $\mat A^{(n)}, n = 1,\cdots, N_1$, along its first $N_1$ modes, combined with a mode-$n$ resource allocation $(n = 1,\cdots, N_1)$ relatively to the mode-$(N_1+1)$ of the transformed tensor $\matc X$, by means of the allocation matrices $\mat \Phi^{(n)}$.
\item In telecommunications applications, the output modes will be called diversity modes because they correspond to time, space and frequency diversities, whereas the input modes are associated with resources like transmit antennas, codes, and data streams. For these applications, the matrices $\mat \Phi^{(n)}$ are formed with 0's and 1's, and they can be interpreted as allocation matrices used for allocating some resources $r_n$ to the output mode-$(N_1+1)$. Another way to take resource allocations into account consists in replacing the $N_1$ allocation matrices $\mat \Phi^{(n)}$  by the $(N_1+1)^{th}$-order allocation tensor $\matc F \inten{R_1\times\cdots}{R_{N_1}}{I_{N_1+1}}$ defined in (\ref{eq:Paratuck allocation tensor}).
\item Special cases:
\begin{itemize}
\item For $N_1 = 2$ and $N = 3$, we obtain the standard PARATUCK-2 model introduced in \cite{Harshman96}. Eq. (\ref{eq:general Paratuck model}) then becomes
\beq \label{eq:Paratuck2 model}
x_{i_1,i_2,i_3}=\sum\limits_{r_1=1}^{R_1}\sum\limits_{r_2=1}^{R_2}c_{r_1,r_2}a^{(1)}_{i_1,r_1}a^{(2)}_{i_2,r_2}\phi^{(1)}_{r_1,i_3}\phi^{(2)}_{r_2,i_3}
\eeq
The allocation tensor $\matc F$ defined in (\ref{eq:Paratuck allocation tensor}) can be rewritten as
\beq \label{eq:PARAFAC allocation}
f_{r_1,r_2,i_3}=\phi^{(1)}_{r_1,i_3}\phi^{(2)}_{r_2,i_3}=\sum\limits_{j=1}^{I_3}\phi^{(1)}_{r_1,j}\phi^{(2)}_{r_2,j}\delta_{i_3,j}
\eeq
which corresponds to a PARAFAC model with matrix factors ($\bf{\Phi}^{(1)},\bf{\Phi}^{(2)},\mat I_{I_3}$). The PARATUCK-2 model (\ref{eq:Paratuck2 model}) can then be viewed as a Tucker-2 model $\matc X=\matc G{\times}_1\mat A^{(1)}{\times}_2\mat A^{(2)}$ with the core tensor $\matc G\inten{R_1}{R_2}{I_3}$ given by the Hadamard product of $\mat C\inmat{R_1}{R_2}$ and $\matc F\inten{R_1}{R_2}{I_3}$ along their common modes $\{r_1,r_2\}$
\beq
\matc G=\mat C\underset{\{r_1,r_2\}}{\odot}\matc F \nonumber
\eeq
This combination of a Tucker-2 model for $\matc X$ with a PARAFAC model for $\matc F$ gave rise to the name PARATUCK-2. The constraint matrices ($\bf{\Phi}^{(1)},\bf{\Phi}^{(2)}$) define interactions between columns of the factor matrices ($\mat A^{(1)},\mat A^{(2)}$), along the mode-3 of $\matc X$, while the matrix $\mat C$ contains the weights of these interactions.
\item For $N_1 = 2$ and $N = 4$, we obtain the PARATUCK-(2,4) model introduced in \cite{FCAR11EUSIPCO}
\beq \label{eq:Paratuck-(2,4) model}
x_{i_1,i_2,i_3,i_4}=\sum\limits_{r_1=1}^{R_1}\sum\limits_{r_2=1}^{R_2}c_{r_1,r_2,i_4}a^{(1)}_{i_1,r_1}a^{(2)}_{i_2,r_2}\phi^{(1)}_{r_1,i_3}\phi^{(2)}_{r_2,i_3}
\eeq
As for the PARATUCK-2 model, the PARATUCK-(2,4) can be viewed as a combination of a Tucker-(2,4) model for $\matc X=\matc G{\times}_1\mat A^{(1)}{\times}_2\mat A^{(2)}\inten{I_1}{I_2}{I_3\times I_4}$ with a core tensor $\matc G\inten{R_1}{R_2}{I_3\times I_4}$ given by the Hadamard product of the tensors $\matc C\inten{R_1}{R_2}{I_4}$ and $\matc F\inten{R_1}{R_2}{I_3}$ along their common modes $\{r_1,r_2\}$
\beq
\matc G=\matc C\underset{\{r_1,r_2\}}{\odot}\matc F \nonumber
\eeq
with the same allocation tensor $\matc F$ defined in (\ref{eq:PARAFAC allocation}).
\end{itemize}
\end{itemize}
\subsection{\bf{Rewriting of PARATUCK models as Constrained PARAFAC Models}}

This rewriting of PARATUCK models as constrained PARAFAC models can be used to deduce both matrix unfoldings by means of the general formula (\ref{eq:Parafac mat rep}), and sufficient conditions for essential uniqueness of such PARATUCK models, as will be shown in Section IV.

\subsubsection{\bf{Link between PARATUCK-(2,4) and constrained PARAFAC-4 models}}

We now establish the link between the PARATUCK-(2,4) model (\ref{eq:Paratuck-(2,4) model})  and the fourth-order constrained PARAFAC model
\beq \label{eq:Paratuck-(2,4) model equiv}
x_{i_1,i_2,i_3,i_4}=\sum\limits_{r=1}^{R}{a}_{i_1,r}{b}_{i_2,r}{f}_{i_3,r}{d}_{i_4,r}\quad \textrm{with}\quad R=R_1R_2
\eeq
whose matrix factors (${\mat A}\inmat {I_1}{R}$, ${\mat B}\inmat {I_2}{R}$, ${\mat F}\inmat {I_3}{R}$, ${\mat D}\inmat {I_4}{R}$), and constraint matrices ($\mat \Psi^{(1)},\mat \Psi^{(2)}$) acting on the original factors ($\mat A^{(1)},\mat A^{(2)}$), are given by
\beq \label{eq:Const Parafac factors}
{\mat A}=\mat A^{(1)} \mat \Psi^{(1)}, \quad {\mat B}=\mat A^{(2)} \mat \Psi^{(2)}, \quad {\mat F}=(\mat \Phi^{(1)} \diamond \mat \Phi^{(2)})^T, \quad {\mat D}=\mat C_{I_4\times R_1R_2}\\
\mat \Psi^{(1)}=\mat I_{R_1}\otimes {\mat 1}^T_{R_2}\inmat{R_1}{R_1R_2} ,\quad \mat \Psi^{(2)}={\mat 1}^T_{R_1}\otimes \mat I_{R_2}\inmat{R_2}{R_1R_2} \label{eq:Paratuck-(2,4) constraints}
\eeq
where $\mat C_{I_4\times R_1R_2}\inmat{I_4}{R_1R_2}$ is a mode-3 unfolded matrix of the tensor $\matc C \inten{R_1}{R_2}{I_4}$.

\begin{proof}
See the Appendix.
\end{proof}

\subsection*{Remarks}

\begin{itemize}
\item Application of the formula (\ref{eq:mode-n unfolding PARAFAC}) to the constrained PARAFAC model (\ref{eq:Paratuck-(2,4) model equiv}), with the matrix factors $({\mat A},{\mat B},{\mat F},{\mat D})=(\mat A^{(1)} \mat \Psi^{(1)}, \mat A^{(2)} \mat \Psi^{(2)}, (\mat \Phi^{(1)} \diamond \mat \Phi^{(2)})^T, \mat C_{I_4\times R_1R_2})$, gives the following flat modes-1 and -2 matrix unfoldings for the PARATUCK-(2,4) model (\ref{eq:Paratuck-(2,4) model})
\beq \label{Constr PARAFAC unfolding a}
\mat X_{I_1\times I_2I_3I_4}&=&\mat A^{(1)} \mat \Psi^{(1)}(\mat A^{(2)} \mat \Psi^{(2)} \diamond \mat F \diamond \mat D)^T\inmat{I_1}{I_2I_3I_4},\nonumber\\
\mat X_{I_2\times I_3I_4I_1}&=&\mat A^{(2)} \mat \Psi^{(2)}(\mat F \diamond \mat D\diamond \mat A^{(1)} \mat \Psi^{(1)})^T\inmat{I_2}{I_3I_4I_1}.\nonumber \label{Constr PARAFAC unfolding b}
\eeq
\item The constrained PARAFAC-4 model (\ref{eq:Paratuck-(2,4) model equiv})-(\ref{eq:Paratuck-(2,4) constraints}) can be written in mode-$n$ products notation as
    \beq \label{eq:Paratuck-(2,4) model equiv c}
    \matc X=\matc I_{4,R}{\times}_1\mat A^{(1)}\mat \Psi^{(1)}{\times}_2\mat A^{(2)}\mat \Psi^{(2)}{\times}_3{\mat F}{\times}_4 {\mat D}.
    \eeq
    Defining the core tensor $\matc G\inten{R_1}{R_2}{I_3\times I_4}$ as
    \beq \label{core tensor}
    \matc G=\matc I_{4,R}{\times}_1\mat \Psi^{(1)}{\times}_2\mat \Psi^{(2)}{\times}_3{\mat F}{\times}_4 {\mat D}
    \eeq
    the constrained PARAFAC-4 model can also be viewed as the following Tucker-(2,4) model
    \beq \label{eq:Paratuck-(2,4) model equiv d}
    \matc X=\matc G{\times}_1\mat A^{(1)}{\times}_2\mat A^{(2)}.
    \eeq
It can also be viewed as a CONFAC-(2,4) model with matrix factors $(\mat A^{(1)}, \mat A^{(2)}, \mat F, \mat D)$, and constraint matrices $\mat \Psi^{(1)}$ and $\mat \Psi^{(2)}$ defined in (\ref{eq:Paratuck-(2,4) constraints}).

\item Choosing $\mathds{S}_1=\{i_1,i_2\}$ and $\mathds{S}_2=\{i_3,i_4\}$, the matrix unfolding (\ref {eq:Parafac mat rep}) of the PARAFAC model (\ref{eq:Paratuck-(2,4) model equiv c}) is given by
\beq \label{unfolding of contracted PARAFAC}
\hspace{-6ex} \mat X_{I_1I_2\times I_3I_4}&=&(\mat A^{(1)} \mat \Psi^{(1)}\diamond \mat A^{(2)} \mat \Psi^{(2)}) \big(\mat F \diamond \mat D\big)^T\nonumber \\
&=&(\mat A^{(1)}\otimes \mat A^{(2)}) \big(\mat F \diamond \mat D\big)^T\inmat{I_1I_2}{I_3I_4}
\eeq
\begin{proof}
{\color{black}
Using the identity (\ref{eq:Kron-Khatri}) gives
\beq
\mat A^{(1)} \mat \Psi^{(1)}\diamond\mat A^{(2)} \mat \Psi^{(2)}=(\mat A^{(1)}\otimes \mat A^{(2)})( \mat \Psi^{(1)} \diamond \mat \Psi^{(2)})
\eeq
Replacing $\mat \Psi^{(1)}$ and $\mat \Psi^{(2)}$ by their expressions (\ref{eq:formula APsi1}) and (\ref{eq:formula APsi2}) leads to
\beq
     \mat \Psi^{(1)} \diamond \mat \Psi^{(2)}&=&(\mat I_{R_1}\otimes {\mat 1}^T_{R_2})\diamond ({\mat 1}^T_{R_1}\otimes \mat I_{R_2})\nonumber \\
     &=&\underset{R_1\,\, \textrm{blocks}}{\underbrace{\left[\begin{array}{ccc} \mat I_{R_2} & & \\
& \ddots &\\
& & \mat I_{R_2} \end{array}\right]}}=\mat I_{R_1R_2} \label{PsiKhatri}
\eeq
which implies
 \beq
 \mat A^{(1)} \mat \Psi^{(1)}\diamond\mat A^{(2)} \mat \Psi^{(2)}=\mat A^{(1)} \otimes \mat A^{(2)}, \label{PsiKhatri}
 \eeq
and consequently Eq. (\ref{unfolding of contracted PARAFAC}) can be deduced.\\
This equation can also be obtained from the equivalent Tucker-(2,4) model (\ref{core tensor})-(\ref{eq:Paratuck-(2,4) model equiv d}) as
\beq \label{Xunfolding}
\mat X_{I_1I_2\times I_3I_4}=(\mat A^{(1)}\otimes \mat A^{(2)})\mat G_{R_1R_2\times I_3I_4}
\eeq
with
\beq
\mat G_{R_1R_2\times I_3I_4}=(\mat \Psi^{(1)} \diamond \mat \Psi^{(2)})(\mat F \diamond \mat D)^T\nonumber
\eeq
Using the identity (\ref{PsiKhatri}), we obtain
\beq \label{coretensorb}
\mat G_{R_1R_2\times I_3I_4}=(\mat F \diamond \mat D)^T
\eeq
and replacing $\mat G_{R_1R_2\times I_3I_4}$ by its expression (\ref{coretensorb}) into  (\ref{Xunfolding}) gives (\ref{unfolding of contracted PARAFAC}).
}

\end{proof}
When the allocation matrices ($\mat \Phi^{(1)}$, $\mat \Phi^{(2)}$) and the input tensor $\matc C$ are known, the matrix factors $(\mat A^{(1)},\mat A^{(2)})$ can be estimated through the LS estimation of their Kronecker product using the matrix unfolding (\ref{unfolding of contracted PARAFAC}).
\item   The product $\phi^{(1)}_{r_1,i_3}\phi^{(2)}_{r_2,i_3}$ in (\ref{eq:Paratuck-(2,4) model}) can be replaced by $f_{i_3,r_1,r_2}$ , which amounts to replace the allocation matrices $\mat \Phi^{(1)}$ and $\mat \Phi^{(2)}$  by the third-order allocation tensor $\matc F\inten{I_3}{R_1}{R_2}$, the matrix ${\mat F}=(\mat \Phi^{(1)} \diamond \mat \Phi^{(2)})^T\inmat{I_3}{R_1R_2}$ being equivalent to $\mat F_{I_3\times R_1R_2}\inmat{I_3}{R_1R_2}$, i.e. a mode-1 flat matrix unfolding of the allocation tensor $\matc F$.
\end{itemize}

\subsubsection{\bf{Link between PARATUCK-2 and constrained PARAFAC-3 models}}
By proceeding in the same way as for the PARATUCK-(2,4) model, it is easy to show that the PARATUCK-2 model (\ref{eq:Paratuck2 model}) is equivalent to a third-order constrained PARAFAC model whose matrix factors ${\mat A}\inmat {I_1}{R}$, ${\mat B}\inmat {I_2}{R}$, and ${\mat F}\inmat {I_3}{R}$, with $R=R_1R_2$, are given by
\beq \label{eq:constr PARAFAC factors}
{\mat A}=\mat A^{(1)} \mat \Psi^{(1)}, \quad {\mat B}=\mat A^{(2)} \mat \Psi^{(2)}, \quad {\mat F}=(\mat \Phi^{(1)} \diamond \mat \Phi^{(2)})^Tdiag(\textrm{vec}(\mat C^T))
\eeq
with the same constraint matrices $\mat \Psi^{(1)}$ and $\mat \Psi^{(2)}$ defined in (\ref{eq:Paratuck-(2,4) constraints}).
By analogy with the PARATUCK-(2,4) model, Eq. (\ref{eq:Paratuck-(2,4) model equiv c}), (\ref{eq:Paratuck-(2,4) model equiv d}), and (\ref{unfolding of contracted PARAFAC}) become for the PARATUCK-2 model
\beq \label{eq:Paratuck-2 model equiv c}
    \matc X&=&\matc I_{3,R}{\times}_1\mat A^{(1)}\mat \Psi^{(1)}{\times}_2\mat A^{(2)}\mat \Psi^{(2)}{\times}_3{\mat F}\nonumber\\
    \label{eq:Paratuck-2 model equiv d}
    &=&\matc G{\times}_1\mat A^{(1)}{\times}_2\mat A^{(2)}
    \eeq
    with the core tensor $\matc G\inten{R_1}{R_2}{I_3}$ defined as
    \beq \label{eq:Paratuck-2 model equiv core tensor}
    \matc G=\matc I_{3,R}{\times}_1\mat \Psi^{(1)}{\times}_2\mat \Psi^{(2)}{\times}_3{\mat F},
    \eeq
    and
    \beq
    \mat X_{I_1I_2\times I_3}=(\mat A^{(1)}\otimes \mat A^{(2)}) {\mat F}^T\inmat{I_1I_2}{I_3}.\nonumber
    \eeq

\subsection*{Remarks}
\begin{itemize}
\item Eq. (\ref{eq:Paratuck-2 model equiv d}) and (\ref{eq:Paratuck-2 model equiv core tensor}) allow interpreting the PARATUCK-2 model as a Tucker-(2,3) model, defined in (\ref{eq:Tucker2})-(\ref{eq:Tucker2 a}).
    If we choose $c_{r_1,r_2}=1, \forall r_k=1,\cdots , R_k$, for $k$=1 and 2, and define the allocation tensor $\matc F\inten{R_1}{R_2}{I_3}$ such as $f_{r_1,r_2,i_3}=\phi^{(1)}_{r_1,i_3}\phi^{(2)}_{r_2,i_3}$, the PARATUCK-2 model (\ref{eq:Paratuck2 model}) becomes the following Tucker-(2,3) model
    \beq \label{eq:Tucker2 model a}
x_{i_1,i_2,i_3}=\sum\limits_{r_1=1}^{R_1}\sum\limits_{r_2=1}^{R_2}f_{r_1,r_2,i_3}a^{(1)}_{i_1,r_1}a^{(2)}_{i_2,r_2}\nonumber
\eeq
and the associated constrained PARAFAC-3 model can be deduced from (\ref{eq:constr PARAFAC factors})
\beq
{\mat A}=\mat A^{(1)} \mat \Psi^{(1)}, \quad {\mat B}=\mat A^{(2)} \mat \Psi^{(2)}, \quad {\mat F}=\mat F_{I_3 \times R_1R_2}=(\mat \Phi^{(1)} \diamond \mat \Phi^{(2)})^T\nonumber
\eeq
with the same constraint matrices $\mat \Psi^{(1)}$ and $\mat \Psi^{(2)}$ as those defined in (\ref{eq:Paratuck-(2,4) constraints}). A block Tucker-(2,3) model transformed into a block constrained PARAFAC-3 model was used in \cite{Almeida07SP} for modeling in an unified way three multiuser wireless communication systems.
    \item Now, we show the equivalence of the expressions (\ref{eq:Paratuck-2 model equiv core tensor}) and (\ref{eq:Hadamard products}) of the core tensor. Applying the formula (\ref{eq:mode-n unfolding PARAFAC}) to the PARAFAC model (\ref{eq:Paratuck-2 model equiv core tensor}) gives
     \beq \label{eq:core unfolding}
     \mat G_{I_3\times R_1R_2}=(\mat \Phi^{(1)} \diamond \mat \Phi^{(2)})^T diag(\textrm{vec}(\mat C^T))(\mat \Psi^{(1)} \diamond \mat \Psi^{(2)})^T.
     \eeq

Using the identity (\ref{PsiKhatri}) in Eq. (\ref{eq:core unfolding}) gives $\mat G_{I_3\times R_1R_2}=(\mat \Phi^{(1)} \diamond \mat \Phi^{(2)})^T diag(\textrm{vec}(\mat C^T))$.
\

    For the formula (\ref{eq:Hadamard products}), with $N=3$ and $N_1=2$, we have
    {\color{black}
    \beq
     \matc G=\matc F\underset{\{r_1,r_2\}}\odot \mat C\nonumber
     \eeq
    or equivalently in terms of matrix Hadamard product
     \beq
     \mat G_{I_3\times R_1R_2}=\mat F_{I_3\times R_1R_2}\odot \mat 1_{I_3}\mat c_{1\times R_1R_2}\nonumber
     \eeq
  }
    with $\mat F_{I_3\times R_1R_2}=(\mat \Phi^{(1)} \diamond \mat \Phi^{(2)})^T$, and $\mat c_{1\times R_1R_2}=\textrm{vec}^T({\mat C}^T)$, which gives
    \beq
     \mat G_{I_3\times R_1R_2}=\mat F_{I_3\times R_1R_2}\odot\left[\begin{array}{c}\textrm{vec}^T({\mat C}^T)\nonumber\\
\vdots \\
\textrm{vec}^T({\mat C}^T) \end{array}\right]\Bigg\}\mat I_3 \, \textrm{rows}
    \eeq


 and consequently $\mat G_{I_3\times R_1R_2}=(\mat \Phi^{(1)} \diamond \mat \Phi^{(2)})^Tdiag(\textrm{vec}(\mat C^T))$, showing the equivalence of the two core tensor expressions (\ref{eq:Paratuck-2 model equiv core tensor}) and (\ref{eq:Hadamard products}).
\end{itemize}


\subsubsection{\bf{Link between PARATUCK-$(N-2,N)$ and constrained PARAFAC-$N$ models}}

Let us consider the PARATUCK-$(N_1,N)$ model (\ref{eq:general Paratuck model}) in the case $N_1=N-2$
\beq \label{eq:general Paratuck model bis}
x_{i_1,\cdots,i_{N_1+1},\cdots,i_N }=\sum\limits_{r_1=1}^{R_1}\cdots \sum\limits_{r_{N_1}=1}^{R_{N_1}}c_{r_1,\cdots,r_{N_1},i_N}\prod\limits_{n=1}^{N_1}a^{(n)}_{i_n,r_n}\phi^{(n)}_{r_n,i_{{N_1}+1}}
\eeq
and let us define the change of variables $r=r_{N_1}+\sum\limits_{n=1}^{{N_1}-1}(r_n - 1)\prod\limits_{i=n+1}^{N_1}R_i$ corresponding to a combination of the $N_1$ modes associated with the constraints/allocations. Eq. (\ref{eq:general Paratuck model bis}) can then be written as the following constrained PARAFAC-$N$ model
\beq \label{eq:gener Paratuck model equiv}
x_{i_1,\cdots,i_N}=\sum\limits_{r=1}^{R} \prod\limits_{n=1}^{N}\bar{a}^{(n)}_{i_n,r},\quad R=\prod\limits_{i=1}^{N_1}R_i
\eeq
with the following matrix factors
\beq \label{eq:equiv gener Paratuck model factors}
\bar{\mat A}^{(n)}=\mat A^{(n)} \mat \Psi^{(n)},\, n=1, \cdots, N_1; \quad {\mat F}={\left(\overset{N}{\underset{n=1}{\diamond}} \mat \Phi^{(n)} \right)}^T; \quad {\mat D}=\mat C_{I_N\times R_1\cdots R_{N_1}},\nonumber
\eeq
where $\mat C_{I_N\times R_1\cdots R_{N_1}}\inmat{I_N}{R_1\cdots R_{N_1}}$ is a mode-$(N_1+1)$ unfolded matrix of the tensor $\matc C\inten {R_1\times\cdots}{R_{N_1}}{I_N}$,
and the constraint matrices are given in (\ref{eq: extension formula}) as
\beq \label{eq:equiv gener Paratuck model constr}
\mat \Psi^{(n)}={\mat 1}^T_{R_1}\otimes\cdots\otimes {\mat 1}^T_{R_{n-1}}\otimes\mat I_{R_n}\otimes{\mat 1}^T_{R_{n+1}}\otimes\cdots \otimes{\mat 1}^T_{R_N}\inmat{R_n}{R}, n= 1, \cdots, N_1.\nonumber
\eeq
\
The constrained PARAFAC model (\ref{eq:gener Paratuck model equiv}) can also be written as a Tucker-$(N_1,N)$ model (\ref{eq:Tucker-N1 a}) with the core tensor defined in (\ref{eq:Hadamard products}), or, equivalently,
\beq
    \matc G=\matc I_{N,R}{\times}_{n=1}^{N-2}\mat \Psi^{(n)}{\times}_{N-1} {\mat F} {\times}_N {\mat D}.\nonumber
    \eeq

\subsection{\bf{Comparison of constrained tensor models}}

To conclude this presentation, we compare the so called CONFAC-$(N_1,N)$ and PARATUCK-$(N_1,N)$ constrained tensor models, introduced in this paper with a resource allocation point of view. Due to the PARAFAC structure (\ref{eq: Paralind core}) of the core tensor of CONFAC models, each element $x_{i_1,\cdots,i_N}$ of the output tensor $\matc X$ is the sum of $R$ components as shown in (\ref{eq: Paralind model c}). Moreover, due to the special structure of the allocation matrices $\mat \Phi^{(n)}$ whose the columns are unit vectors, each component $r$ is the result of a combination of $N$ resources, under the form of the product $\prod \limits_{n=1}^{N}a^{(n)}_{i_n,r_n}$, the $N$ resources being fixed by the allocation matrices $\mat \Phi^{(n)}\inmat{R_n}{R}$.


With the CONFAC-$(N_1,N)$ model (\ref {eq: Paralind N_1 model a}), each component $r$ is a combination of $N_1$ resources $(r_1,\cdots ,r_{N_1})$ determined by the allocation matrices $\mat \Phi^{(n)}\inmat{R_n}{R}$ for $n=1,\cdots ,N_1$.


There are two main differences between the PARATUCK-$(N_1,N)$ models (\ref{eq:general Paratuck model}) and the CONFAC models (\ref {eq: Paralind model a}). The first one is that the allocation matrices of PARATUCK models, formed with 0's and 1's, have not necessarily unit vectors as column vectors, which means that it is possible to allocate $\gamma_n=\sum\limits_{r_n=1}^{R_n}\phi^{(n)}_{r_n,i_{N_1+1}}$ resources $r_n$ to the $(N_1+1)^{th}$-mode of the output tensor $\matc X$. The second one results from the interpretation of PARATUCK-$(N_1,N)$ models as Tucker-$(N_1,N)$ models, implying that each element $x_{i_1,\cdots,i_N}$ of $\matc X$ is equal to the sum of $\sum\limits_{r_1=1}^{R_1}\cdots \sum\limits_{r_{N_1}=1}^{R_{N_1}}f_{r_1,\cdots,r_{N_1},i_{{N_1}+1}}$ terms, where $f_{r_1,\cdots,r_{N_1},i_{{N_1}+1}}$ is an entry of the allocation tensor $\matc F$ defined in (\ref{eq:Paratuck allocation tensor}), each term being a combination of resources under the form of products $\prod \limits_{n=1}^{N_1}a^{(n)}_{i_n,r_n}$. Moreover, in telecommunication applications, the input tensor $\matc C$ can be used as a code tensor.\\

{\color{black}Another way to compare PARALIND/CONFAC and PARATUCK models is in terms of dependencies/interactions between their factor matrices. In the case of PARALIND/CONFAC models, as pointed out by Eq. (\ref{eq: Paralind model a}), the constraint matrices act independently on each factor matrix, expliciting linear dependencies between columns of these matrices. For PARATUCK models, their writing as Tucker-$(N_1,N)$ models with the core tensor defined in (\ref{eq:Hadamard products}) allows to interpret the tensor $\matc F$ as an interaction tensor which defines interactions between $N_1$ factor matrices, the tensor $\matc C$ providing the strength of these interactions.
}
\\

The main constrained PARAFAC models are summarized in Tables I and II.

\begin{table}[!h]
{\scriptsize
\begin{center} \caption{Main tensor models}
\begin{tabular}{|c|c|c|} 
\hline
 Models & Scalar writings &  mode-$n$ product based writings \\
 \hline \hline
  PARAFAC-3 &$x_{i_1,i_2,i_3}=\sum\limits_{r}^{R}a^{(1)}_{i_1,r}a^{(2)}_{i_2,r}a^{(3)}_{i_3,r}$ &  $\matc X=\matc I_{3,R} \times_1\mat A^{(1)}\times_2\mat A^{(2)}\times_3\mat A^{(3)}$  \\
   & &  \\
  \hline
 Tucker-3 &
 $x_{i_1,i_2,i_3}=\sum\limits_{r_1=1}^{R_1}\sum\limits_{r_2=1}^{R_2}\sum\limits_{r_3=1}^{R_3}g_{r_1,r_2,r_3}a^{(1)}_{i_1,r_1}a^{(2)}_{i_2,r_2}a^{(3)}_{i_3,r_3}$ & $\matc X=\matc G \times_1\mat A^{(1)}\times_2\mat A^{(2)}\times_3\mat A^{(3)}$  \\
   & &  \\
  \hline
 Tucker-(2,3) & $x_{i_1,i_2,i_3}=\sum\limits_{r_1=1}^{R_1}\sum\limits_{r_2=1}^{R_2}g_{r_1,r_2,i_3}a^{(1)}_{i_1,r_1}a^{(2)}_{i_2,r_2}$ & $\matc X=\matc G \times_1\mat A^{(1)}\times_2\mat A^{(2)}$  \\
   & &  \\
  \hline
 PARALIND/ & $x_{i_1,i_2,i_3}=\sum\limits_{r_1=1}^{R_1}\sum\limits_{r_2=1}^{R_2}\sum\limits_{r_3=1}^{R_3}g_{r_1,r_2,r_3}a^{(1)}_{i_1,r_1}a^{(2)}_{i_2,r_2}a^{(3)}_{i_3,r_3}$ & $\matc X=\matc G \times_1\mat A^{(1)}\times_2\mat A^{(2)}\times_3\mat A^{(3)}$ \\
 CONFAC-3  & & \\
  & $g_{r_1,r_2,r_3}=\sum\limits_{r}^{R}\varphi^{(1)}_{r_1,r}\varphi^{(2)}_{r_2,r}\varphi^{(3)}_{r_3,r}$ & $\matc G=\matc I_{3,R} \times_1\mat \Phi^{(1)}\times_2\mat \Phi^{(2)}\times_3\mat \Phi^{(3)}$ \\
    & &   \\
 \hline
 Paratuck-2 & $x_{i_1,i_2,i_3}=\sum\limits_{r_1=1}^{R_1}\sum\limits_{r_2=1}^{R_2}g_{r_1,r_2,i_3}a^{(1)}_{i_1,r_1}a^{(2)}_{i_2,r_2}$ & $\matc X=\matc G \times_1\mat A^{(1)}\times_2\mat A^{(2)}$    \\
 & $g_{r_1,r_2,i_3}=c_{r_1,r_2}\varphi^{(1)}_{r_1,i_3}\varphi^{(2)}_{r_2,i_3}$  & $\matc G=\matc I_{3,R}{\times}_1\mat \Psi^{(1)}{\times}_2\mat \Psi^{(2)}{\times}_3\bar{\mat C}$ \quad , \quad $R=R_1R_2$\\
   & & $\mat \Psi^{(1)}=\mat I_{R_1}\otimes {\mat 1}^T_{R_2} \quad,\quad \mat \Psi^{(2)}={\mat 1}^T_{R_1}\otimes \mat I_{R_2}$ \\
   & & $\bar{\mat C}=(\mat \Phi^{(1)} \diamond \mat \Phi^{(2)})^Tdiag(\textrm{vec}(\mat C^T))$ \\
\hline
Paratuck-(2,4) & $x_{i_1,i_2,i_3,i_4}=\sum\limits_{r_1=1}^{R_1}\sum\limits_{r_2=1}^{R_2}g_{r_1,r_2,i_3,i_4}a^{(1)}_{i_1,r_1}a^{(2)}_{i_2,r_2}$ & $\matc X=\matc G \times_1\mat A^{(1)}\times_2\mat A^{(2)}$    \\
 & $g_{r_1,r_2,i_3,i_4}=c_{r_1,r_2,i_4}\varphi^{(1)}_{r_1,i_3}\varphi^{(2)}_{r_2,i_3}$  &  $\matc G=\matc I_{4,R}{\times}_1\mat \Psi^{(1)}{\times}_2\mat \Psi^{(2)}{\times}_3\bar{\mat F}{\times}_4 \bar{\mat D}$ \quad , \quad $R=R_1R_2$ \\
  & & $\mat \Psi^{(1)}=\mat I_{R_1}\otimes {\mat 1}^T_{R_2} \quad,\quad \mat \Psi^{(2)}={\mat 1}^T_{R_1}\otimes \mat I_{R_2}$ \\
   & & $\bar{\mat F}=(\mat \Phi^{(1)} \diamond \mat \Phi^{(2)})^T \quad,\quad  \bar{\mat D}=\mat C_{I_4\times R_1R_2}$ \\
 \hline
\end{tabular}
\end{center}
}
\end{table}

\begin{table}[!h]
{\scriptsize
\begin{center} \caption{Equivalent constrained PARAFAC models}
\begin{tabular}{|c|c|c|} 
\hline
 Models & Equivalent constrained PARAFAC model &  Matrix unfoldings \\
 \hline \hline 

 & &   \\  PARAFAC-3 
   & &  $\mat X_{I_1\times I_2I_3}=\mat A^{(1)}(\mat A^{(2)}\diamond\mat A^{(3)})^T$\\
   & &  \\

  \hline 
& &   \\  Tucker-3 
   &  & $\mat X_{I_1\times I_2I_3}=\mat A^{(1)}\mat G_{R_1\times R_2R_3}(\mat A^{(2)}\otimes\mat A^{(3)})^T$ \\
   & &  \\
  \hline 
 & &   \\ Tucker-(2,3) 
 & $\matc X=\matc I_{3,R} \times_1\mat A^{(1)}\mat \Psi^{(1)}\times_2\mat A^{(2)}\mat \Psi^{(2)}\times_3\mat G_{I_3\times R_1R_2}$ & $\mat X_{I_1\times I_2I_3}=\mat A^{(1)}\mat \Psi^{(1)}(\mat A^{(2)}\mat \Psi^{(2)}\diamond\mat G_{I_3\times R_1R_2})^T$ \\
   &  $\mat \Psi^{(1)}=\mat I_{R_1}\otimes {\mat 1}^T_{R_2} \quad,\quad \mat \Psi^{(2)}={\mat 1}^T_{R_1}\otimes \mat I_{R_2}$ & \\
   & &  \\
  \hline 
    & &   \\
 PARALIND/ & $\matc X=\matc I_{3,R} \times_1\mat A^{(1)}\mat \Phi^{(1)}\times_2\mat A^{(2)}\mat \Phi^{(2)}\times_3\mat A^{(3)}\mat \Phi^{(3)}$ & $\mat X_{I_1\times I_2I_3}=\mat A^{(1)}\mat \Phi^{(1)}(\mat A^{(2)}\mat \Phi^{(2)}\diamond\mat A^{(3)}\mat \Phi^{(3)})^T$  \\
  CONFAC-3
  & & =$\mat A^{(1)}\mat \Phi^{(1)}(\mat \Phi^{(2)}\diamond\mat \Phi^{(3)})^T(\mat A^{(2)}\otimes\mat A^{(3)})^T$  \\
    & &  \\
 \hline 
& &   \\  Paratuck-2 & $\matc X=\matc I_{3,R}{\times}_1\mat A^{(1)}\mat \Psi^{(1)}{\times}_2\mat A^{(2)}\mat \Psi^{(2)}{\times}_3\bar{\mat C}$ &  \\
  & $\mat \Psi^{(1)}=\mat I_{R_1}\otimes {\mat 1}^T_{R_2} \quad,\quad \mat \Psi^{(2)}={\mat 1}^T_{R_1}\otimes \mat I_{R_2}$ & $\mat X_{I_1\times I_2I_3}=\mat A^{(1)}\mat \Psi^{(1)}(\mat A^{(2)}\mat \Psi^{(2)}\diamond \bar{\mat C})^T$   \\
   &  $\bar{\mat C}=(\mat \Phi^{(1)} \diamond \mat \Phi^{(2)})^Tdiag(\textrm{vec}(\mat C^T))$ & \\
 & & \\
\hline 
& &   \\ Paratuck-(2,4) & $\matc X=\matc I_{4,R}{\times}_1\mat A^{(1)}\mat \Psi^{(1)}{\times}_2\mat A^{(2)}\mat \Psi^{(2)}{\times}_3\bar{\mat F}{\times}_4 \bar{\mat D}$ & $\mat X_{I_1\times I_2I_3I_4}=\mat A^{(1)} \mat \Psi^{(1)}(\mat A^{(2)} \mat \Psi^{(2)} \diamond \bar{\mat F} \diamond \mat C_{I_4\times R_1R_2})^T$  \\
  &  $\mat \Psi^{(1)}=\mat I_{R_1}\otimes {\mat 1}^T_{R_2} \quad,\quad \mat \Psi^{(2)}={\mat 1}^T_{R_1}\otimes \mat I_{R_2}$ & $\mat X_{I_1I_2\times I_3I_4}=(\mat A^{(1)}\otimes \mat A^{(2)})$  \\
 & $\bar{\mat F}=(\mat \Phi^{(1)} \diamond \mat \Phi^{(2)})^T \quad,\quad  \bar{\mat D}=\mat C_{I_4\times R_1R_2}$ & $\big((\mat \Phi^{(1)} \diamond \mat \Phi^{(2)})^T \diamond \mat C_{I_4\times R_1R_2}\big)^T$ \\
  & & \\
 \hline
\end{tabular}
\end{center}
}
\end{table}

\section{\bf{Uniqueness Issue}}

{\color{black}Several results exist for essential uniqueness of PARAFAC models, i.e. uniqueness of factor matrices up to column permutation and scaling. These results concern both deterministic and generic uniqueness, i.e. uniqueness for a particular PARAFAC model, or uniqueness with probability one in the case where the entries of the factor matrices are drawn from continuous distributions. An overview of main uniqueness conditions of PARAFAC models of third-order tensors can be found in \cite{DL2013} for the deterministic case, and in \cite{DL2014} for the generic case. Hereafter, we briefly summarized some basic results on uniqueness of PARAFAC models. The case with linearly dependent loadings is also discussed.} Then, we present new results concerning the uniqueness of PARATUCK models. These results are directly deduced from sufficient conditions for essential uniqueness of their associated constrained PARAFAC models, as established in the previous section. These conditions involving the notion of $k$-rank of a matrix, we first recall the definition of $k$-rank.

\

\emph{Definition of $k$-rank}
\

The $k$-rank (also called Kruskal's rank) of a matrix $\mat A\inmat{I}{R}$, denoted by $k_{\mat A}$, is the largest integer such that any set of $k_{\mat A}$ columns of $\mat A$ is linearly independent.\\
{\color{black}It is obvious that $k_{\mat A}\leq r_{\mat A}$.
}

\subsection{\bf{Uniqueness of PARAFAC-$N$ models \cite{Sid00Nway}}}

The PARAFAC-$N$ model (\ref{eq: Parafac model a})-(\ref{eq: Parafac model c}) is essentially unique, i.e. its factor matrices $\mat A^{(n)}\inmat{I_n}{R}, n= 1,\cdots,N$, are unique up to column permutation and scaling, if
\beq \label{eq:gener Kruskal cond}
 \sum\limits_{n=1}^N k_{\mat A^{(n)}}\geq 2R+N-1
\eeq
\
{\color{black}
Essential uniqueness means that two sets of factor matrices are linked by the following relations ${\hat{\mat A}}^{(n)}=\mat A^{(n)}\bf{\Pi} \bf{\Lambda}^{(n)}$, for $n=1,\cdots,N$, where $\bf{\Pi}$ is a permutation matrix, and $\bf{\Lambda}^{(n)}$ are nonsingular diagonal matrices such as $\prod\limits_{n=1}^N \bf{\Lambda}^{(n)}=\mat I_R$.\\
}
{\color{black}In the generic case, the factor matrices are full rank, which implies $k_{\mat A^{(n)}}=\textrm{min}(I_n,R)$, and the Kruskal's condition (\ref{eq:gener Kruskal cond}) becomes
\beq
\sum\limits_{n=1}^N \textrm{min}(I_n,R)\geq 2R+N-1
\eeq
}

\

\emph{\bf{Case of third-order PARAFAC models}}
\

Consider a third-order tensor $\matc X\inten{I}{J}{K}$ of rank $R$, satisfying a PARAFAC model with matrix factors $(\mat A,\mat B,\mat C)$. The Kruskal's condition (\ref{eq:gener Kruskal cond}) becomes
\beq \label{eq:Kruskal cond}
k_{\mat A}+k_{\mat B}+k_{\mat C}\geq 2R+2
\eeq

\ 				 	
\textbf{Remarks}
\

\begin{itemize}
\item The condition (\ref{eq:gener Kruskal cond}) is sufficient but not necessary for essential uniqueness. This condition does not hold when $R = 1$. It is also necessary for $R = 2$ and $R = 3$ but not for $R > 3$. See \cite{tenBerge022}.\
\item The first sufficient condition for essential uniqueness of third-order PARAFAC models was established by Harshman in \cite{Harshman72}, then generalized by Kruskal in \cite{Kru77} using the concept of $k$-rank. A more accessible proof of Kruskal's condition is provided in \cite{Stegeman07}. The Kruskal's condition was extended to complex-valued tensors in \cite{Sid00} and to $N$-way arrays, with $N > 3$, in \cite{Sid00Nway}.\
\item Necessary and sufficient uniqueness conditions more relaxed than the Kruskal's one were established for third- and fourth-order tensors, under the assumption that at least one matrix factor is full column-rank \cite{JiangSid04}, \cite{DeLathauwer2006}. These conditions are complicated to apply. Other more relaxed conditions have been recently derived, independently by Stegeman \cite{Stegeman08} and Guo et al. \cite{GMB11}, for third-order PARAFAC models with a full column-rank matrix factor.\
\item From the condition (\ref{eq:Kruskal cond}), we can conclude that, if two matrix factors ($\mat A$ and $\mat B$) are full column rank ($k_A=k_B=R$) , then the PARAFAC model is essentially unique if the third matrix factor ($\mat C$) has no proportional columns ($k_C>1$).\
\item If one matrix factor ($\mat C$ for instance) is full column rank, then (\ref{eq:Kruskal cond}) gives
\be \label{eq:Kruskal cond a}
k_{\mat A}+k_{\mat B}\geq R+2
\ee
\
In \cite{Stegeman08} and \cite{GMB11}, it is shown that the PARAFAC model ($\mat A,\mat B,\mat C$), with $\mat C$ of full column rank, is essentially unique if the other two matrix factors $\mat A$ and $\mat B$ satisfy the following conditions
\beq
&&1) \quad k_{\mat A},k_{\mat B}\geq 2 \nonumber \\
\label{eq:Kruskal cond b}
&&2) \quad r_{\mat A}+k_{\mat B}\geq R+2\quad \textrm{or}\quad r_{\mat B}+k_{\mat A}\geq R+2
\eeq
\
Conditions (\ref{eq:Kruskal cond b}) are more relaxed than (\ref{eq:Kruskal cond a}). Indeed, if for instance $k_{\mat A}=2$ and $r_{\mat A}=k_{\mat A}+\delta$ with $\delta>0$, application of (\ref{eq:Kruskal cond a}) implies $k_{\mat B}=R$, i.e. $\mat B$ must be full column rank, whereas (\ref{eq:Kruskal cond b}) gives $k_{\mat B}\geq R-\delta$ which does not require that $\mat B$ be full column rank.\\
\item {\color{black}When one matrix factor ($\mat C$ for instance) is known and the Kruskal's condition (\ref{eq:Kruskal cond}) is satisfied, as it is often the case in telecommunication applications, essential uniqueness is ensured without permutation ambiguity and with only scaling ambiguities ($\bf{\Lambda}_{\mat A},\bf{\Lambda}_{\mat B}$) such as $\bf{\Lambda}_{\mat A}\bf{\Lambda}_{\mat B}=\mat I_R$.}
\end{itemize}

\subsection{Uniqueness of PARAFAC models with linearly dependent loadings}

If one matrix factor contains at least two proportional columns, i.e. its $k$-rank is equal to one, then the Kruskal's condition (\ref{eq:Kruskal cond}) cannot be satisfied. In this case, partial uniqueness can be ensured, i.e. some columns of some matrix factors are essentially unique while the others are unique up to multiplication by a non-singular matrix \cite{B2004}. To illustrate this result, let us consider the case of the PARAFAC model of a fourth-order tensor $\matc X\inten{I\times J}{K}{L}$ with factor matrices $(\mat A,\mat B,\mat C,\mat D)$ whose two of them have two identical columns, at the same position
    \beq
    \mat A=\left[\begin{array}{ccc}\mat A_1\ & \mat a \ & \mat a \end{array}\right], \mat B=\left[\begin{array}{ccc}\mat B_1\ & \mat b & \mat b \end{array}\right], \mat C=\left[\begin{array}{cc}\mat C_1\ & \mat C_2 \end{array}\right], \mat D=\left[\begin{array}{cc}\mat D_1\ & \mat D_2 \end{array}\right]\nonumber
    \eeq
    \
    with $\mat A_1\inmat {I}{(R-2)}, \mat a\inmat {I}{1}, \mat B_1\inmat {J}{(R-2)}, \mat b\inmat {J}{1}, \mat C_1\inmat {K}{(R-2)}, \mat C_2\inmat {K}{2}, \mat D_1\inmat {L}{(R-2)},\mat D_2\inmat {L}{2}$. {\color{black} We have $k_{\mat A}=k_{\mat B}=1$, and consequently the uniqueness condition (\ref{eq:gener Kruskal cond}) for $N$=4 becomes $k_{\mat C}+k_{\mat D}\geq 2R+1$, which cannot be satisfied. In this case, we have partial uniqueness.} Indeed, the matrix slices (\ref{eq: fourth order slices}) can be developed as follows
    \beq
    \mat X_{ij..}&=&\mat C D_j(\mat B) D_i(\mat A)\mat D^T=\left[\begin{array}{cc}\mat C_1 & \mat C_2 \end{array}\right] \left[\begin{array}{cc}
    D_j(\mat B_1) D_i(\mat A_1) & \mat 0_{(R-2)\times{2}}\\
    \mat 0_{{2}\times (R-2)} & a_i b_j\mat I_2 \end{array}\right]\left[\begin{array}{c}\mat D_1^T\\
    \mat D_2^T \end{array}\right]\nonumber\\
    &=&\mat C_1 D_j(\mat B_1) D_i(\mat A_1)\mat D_1^T + a_i b_j \mat C_2 \mat D_2^T.\nonumber
    \eeq
    From this expression, it is easy to conclude that the last two columns of $\mat C$ and $\mat D$ are unique up to a rotational indeterminacy. Indeed, if one replaces the matrices $(\mat C_2,\mat D_2)$ by $(\mat C_2\mat T,\mat D_2\mat T^{-T})$, where $\mat T\inmat{2}{2}$ is a non-singular matrix, the matrix slices $\mat X_{ij..}$ remain unchanged. So, the PARAFAC model is said partially unique in the sense that only the blocks ($\mat A_1, \mat B_1, \mat C_1, \mat D_1$) are essentially unique, the blocks $\mat C_2$  and  $\mat D_2$ being unique up to a non-singular matrix. Essential uniqueness means that any alternative blocks ($\hat{\mat A}_1, \hat{\mat B}_1, \hat{\mat C}_1, \hat{\mat D}_1$) are such as $\hat{\mat A}_1=\mat A_1 \mat \Pi \mat \Delta_a, \hat{\mat B}_1=\mat B_1 \mat \Pi \mat \Delta_b, \hat{\mat C}_1=\mat C_1 \mat \Pi \mat \Delta_c, \hat{\mat D}_1=\mat D_1 \mat \Pi \mat \Delta_d$, where $\mat \Pi$ is a permutation matrix, and $\mat \Delta_a$, $\mat \Delta_b$, $\mat \Delta_c$, and $\mat \Delta_d$ are diagonal matrices such as $\mat \Delta_a\mat \Delta_b\mat \Delta_c\mat \Delta_d = \mat I_{R-2}$. {\color{black} In \cite{BMC2011}, sufficient conditions are provided for essential uniqueness of fourth-order CP models with one full column rank factor matrix, and at most three collinear factor matrices, i.e. having one (or more) column(s) proportional to another column. Uniqueness is ensured if any pair of proportional columns can not be common to two collinear factors, which is not the case of the example above due to the fact that the last two columns of $\mat A$ and $\mat B$ are assumed to be equal. }
    \

{\color{black}

The PARALIND and CONFAC models represent a class of constrained PARAFAC models where the columns of one or more matrix factors are linearly dependent or collinear. In the case of CONFAC models, such a collinearity takes the form of repeated columns that are explicitly modeled by means of constraint matrices. The work \cite{StegAlm09} derived both essential uniqueness conditions and partial uniqueness conditions for PARALIND/CONFAC models of third-order tensors. Therein, the relation with uniqueness of constrained Tucker3 models and the block decomposition in rank-($L$,$L$,1) terms is also discussed. The essential uniqueness condition for a given matrix factor in PARALIND models makes use of Kruskal's Permutation Lemma \cite{Kru77,JiangSid04}.

Consider a third-order tensor $\matc X\inten{I}{J}{K}$ satisfying a PARALIND model with matrix factors $(\mat A,\mat B,\mat C)$, and constraint matrices $\bm{\Phi}^{(i)}$, $i=1,2,3$. 
Suppose $(\mat B\otimes \mat C)\mat G_{R_2R_3 \times R_1}$ and $\mat A$ have full column rank and let $\omega(\cdot)$ denote the number of nonzero elements of its vector argument. Define $N_i= rank{(\bm{\Phi}^{(2)}\,diag(\bm{\Phi}^{(1)}_{i,\,.})\,\bm{\Phi}^{(3)T})}$, $i=1,\ldots, R_1$. If for any vector $\mat d$,
\beq \label{eq-cond}
&&rank{\left[\mat B\bm{\Phi}^{(2)}\,diag(\mat d^T\bm{\Phi}^{(1)})\,(\mat C\bm{\Phi}^{(3)})^T\right]}\leq max(N_1,\ldots, N_{R_1})\nonumber\\
&&\nonumber\\
&& \qquad\qquad \qquad \qquad \quad\quad {\rm implies} \quad \omega(\mat d)\leq 1\,
\eeq
then $\mat A$ is essentially unique \cite{StegAlm09}.
The uniqueness condition for $\mat B$ and $\mat C$ is analogous to condition (\ref{eq-cond}) by interchanging the roles of $\bm{\Phi}^{(1)}$, $\bm{\Phi}^{(2)}$ and $\bm{\Phi}^{(3)}$.

When PARALIND model reduces to PARAFAC model, condition (\ref{eq-cond}) is identical to Condition B of \cite{JiangSid04} for the essential uniqueness of the PARAFAC model in the case of a full column rank matrix factor. More recently in \cite{StegLam12}, improved versions of the main uniqueness conditions of PARALIND/CONFAC models have been derived. The results presented therein involve simpler proofs than those of \cite{StegAlm09}. Moreover, the associated  uniqueness conditions are easy-to-check in comparison with the ones presented earlier in \cite{StegAlm09}.

In \cite{GuoMironBrieSteg12}, a ``uni-mode'' uniqueness condition is derived for a PARAFAC model with linearly dependent (proportional/identical) columns in one matrix factor. This condition is particularly useful for a subclass of PARALIND/CONFAC models with $\bm{\Phi}^{(2)}=\bm{\Phi}^{(3)}=\mat I_{R}$, i.e. when collinearity is confined within the first matrix factor. Let $\bar{\mat A}=\mat A\bm{\Phi}^{(1)}$, where $\bar{\mat A} \inmat{I_1}{R}$ contains collinear columns, the collinearity pattern being captured by $\bm{\Phi}^{(1)}$. Assuming that $\bar{\mat A}$ does not contain an all-zero column, if
\be
r_{\bar{\mat A}} + k_{\mat B} + k_{\mat C} \geq 2R+2,
\ee
then $\bar{\mat A}$ is essentially unique \cite{GuoMironBrieSteg12}. Generalizations of this condition can be obtained by imposing additional constraints on the ranks and $k$-ranks of the matrix factors (see \cite{GuoMironBrieSteg12} for details).

In \cite{BMC2011}, the attention is drawn to the case of fourth-order PARAFAC models with collinear loadings in at most three modes. Note that this type of model can be interpreted as a fourth-order CONFAC model with constraints on the first, second, and third matrix factors. Although collinearity is not explicitly modeled by means of constraint matrices, the uniqueness result of \cite{BMC2011} directly apply to fourth-order CONFAC models.

}

\
\subsection{\bf{Uniqueness of Tucker models}}
\
Contrary to PARAFAC models, the Tucker ones are generally not essentially unique. Indeed, the parameters of Tucker models can be only estimated up to nonsingular transformations characterized by nonsingular matrices $\mat T^{(n)}$ that act on the mode-$n$ matrix factors $\mat A^{(n)}$, and can be cancelled in replacing the core tensor by $\matc G{{\times}^{N}_{n=1}}[\mat T^{(n)}]^{-1}$. This result is easy to verify by applying the property (\ref{eq: mode-n prod property}) of mode-{n} product
\beq
\matc G{{\times}^{N}_{n=1}}[\mat T^{(n)}]^{-1}{{\times}^{N}_{n=1}}\mat A^{(n)}\mat T^{(n)}&=&\matc G{{\times}^{N}_{n=1}}\mat A^{(n)}\mat T^{(n)}[\mat T^{(n)}]^{-1}\nonumber \\
&=&\matc G{{\times}^{N}_{n=1}}\mat A^{(n)}.\nonumber
\eeq

Uniqueness can be obtained by imposing some constraints on the core tensor or the matrix factors. See \cite{Smilde_book} for a review of main results concerning uniqueness of Tucker models, with discussion of three different approaches for simplifying core tensors so that uniqueness is ensured. Uniqueness can also result from a core with information redundancy and structure constraints as in \cite{FBA12} where the core is characterized by matrix slices in Hankel and Vandermonde forms.

\subsection{\bf{Uniqueness of the PARATUCK-(2,4) model}}

Let us consider the PARATUCK-(2,4) model defined by Eq. (\ref{eq:Paratuck-(2,4) model}), with matrix factors $\mat A^{(1)}$ and $\mat A^{(2)}$, constraint matrices $\mat \Phi^{(1)}$ and $\mat \Phi^{(2)}$, and core tensor $\matc C$. As previously shown, this model is equivalent to the constrained PARAFAC model (\ref{eq:Paratuck-(2,4) model equiv}) whose matrix factors are
\beq
{\mat A}=\mat A^{(1)} \mat \Psi^{(1)}, \quad {\mat B}=\mat A^{(2)} \mat \Psi^{(2)}, \quad {\mat F}=(\mat \Phi^{(1)} \diamond \mat \Phi^{(2)})^T, \quad {\mat D}=\mat C_{I_4\times R_1R_2}\nonumber
\eeq
\
with $\mat \Psi^{(1)}$ and $\mat \Psi^{(2)}$ defined in (\ref{eq:Paratuck-(2,4) constraints}). Due to the repetition of some columns of $\mat A^{(1)}$ and $\mat A^{(2)}$, and assuming that these matrices do not contain an all-zero column, we have $k_{\mat A}=k_{\mat B}=1$, and application of the Kruskal's condition (\ref{eq:gener Kruskal cond}), with $N = 4$, gives
\beq
k_{{\mat A}}+k_{{\mat B}}+k_{{\mat F}}+k_{{\mat D}}\geq 2R_1R_2+3\quad \Rightarrow\quad k_{{\mat F}}+k_{{\mat D}}\geq 2R_1R_2+1,\nonumber
\eeq
which can never be satisfied. However, more relaxed sufficient conditions can be established for essential uniqueness of the PARATUCK-(2,4) model. For that purpose, we consider the contracted constrained PARAFAC model obtained by combining the first two modes and using (\ref{PsiKhatri}), which leads to a third-order PARAFAC model with matrix factors
\beq \label{eq:contracted constr Paratuck(2,4)}
({\mat A}\diamond{\mat B},{\mat F},{\mat D})=(\mat A^{(1)} \otimes\mat A^{(2)} ,(\mat \Phi^{(1)} \diamond \mat \Phi^{(2)})^T,\mat C_{I_4\times R_1R_2})
\eeq
{\color{black}
Note that uniqueness of the matrix factors of the contracted PARAFAC model (\ref{eq:contracted constr Paratuck(2,4)}) implies the uniqueness of the matrix factors $\mat A^{(1)}$ and $\mat A^{(2)}$ of the original PARATUCK-(2,4) model. This comes from the fact that $\mat A^{(1)}$ and $\mat A^{(2)}$ can be recovered (up to a scaling factor) from their Kronecker product \cite{VanLoanPitsianis}. Application of the conditions (\ref{eq:Kruskal cond b}) to the contracted PARAFAC model (\ref{eq:contracted constr Paratuck(2,4)}) allows deriving the following theorem.}

\
\emph{\bf{Theorem:}}
\

The PARATUCK-(2,4) model defined by Eq. (\ref{eq:Paratuck-(2,4) model}) is essentially unique
\begin{itemize}
\item 1) When $\mat A^{(1)}$ and $\mat A^{(2)}$ are full column-rank ($r_{\mat A^{(1)}\otimes\mat A^{(2)}}=R_1R_2\Rightarrow k_{\mat A^{(1)}\otimes\mat A^{(2)}}=R_1R_2$)\\
    If $\left\{\begin{array}{ll} k_{(\bm{\Phi}^{(1)}\diamond \bm{\Phi}^{(2)})^T}\geq 2\\
k_{\mat C_{I_4\times R_1R_2}}\geq 2 \end{array}\right.$ and $\left\{\begin{array}{lll} r_{(\bm{\Phi}^{(1)}\diamond \bm{\Phi}^{(2)})^T}+k_{\mat C_{I_4\times R_1R_2}}\geq R_1R_2+2\\
{\rm or}\\
r_{\mat C_{I_4\times R_1R_2}}+k_{(\bm{\Phi}^{(1)}\diamond \bm{\Phi}^{(2)})^T}\geq R_1R_2+2
 \end{array}\right.$\\
\item 2) When $(\bm{\Phi}^{(1)} \diamond \bm{\Phi}^{(2)})^T$ is full column-rank\\
If $\left\{\begin{array}{ll} k_{\mat A^{(1)}\otimes\mat A^{(2)}}\geq 2\\
k_{\mat C_{I_4\times R_1R_2}}\geq 2 \end{array}\right.$ and $\left\{\begin{array}{lll} r_{\mat A^{(1)}}r_{\mat A^{(2)}}+k_{\mat C_{I_4\times R_1R_2}}\geq R_1R_2+2\\
{\rm or}\\
r_{\mat C_{I_4\times R_1R_2}}+k_{\mat A^{(1)}\otimes\mat A^{(2)}}\geq R_1R_2+2
 \end{array}\right.$\\
\item 3) When $\mat C_{I_4\times R_1R_2}$ is full column-rank\\
If $\left\{\begin{array}{ll} k_{\mat A^{(1)}\otimes\mat A^{(2)}}\geq 2\\
k_{(\bm{\Phi}^{(1)}\diamond \bm{\Phi}^{(2)})^T}\geq 2 \end{array}\right.$ and $\left\{\begin{array}{lll} r_{\mat A^{(1)}}r_{\mat A^{(2)}}+k_{(\bm{\Phi}^{(1)}\diamond \bm{\Phi}^{(2)})^T}\geq R_1R_2+2\\
{\rm or}\\
r_{(\bm{\Phi}^{(1)}\diamond \bm{\Phi}^{(2)})^T}+k_{\mat A^{(1)}\otimes\mat A^{(2)}}\geq R_1R_2+2
\end{array}\right.$
\end{itemize}
\

{\color{black}
In \cite{FCAR12SP}, an application of the PARATUCK-(2,4) model to tensor space-time (TST) coding is considered. Therein, the matrix factors $\mat A^{(1)}$ and $\mat A^{(2)}$ represent the symbol and channel matrices to be estimated while the constraint matrices $\bm{\Phi}^{(1)}$ and $\bm{\Phi}^{(2)}$ play the role of allocation matrices of the transmission system and the tensor $\matc C$ is the coding tensor. In this context, $\bm{\Phi}^{(1)}$, $\bm{\Phi}^{(2)}$ and $\matc C$ can be properly designed to satisfy the sufficient conditions of  item 1) of the Theorem. }\\
 The sufficient conditions of this Theorem can easily be extended to the case of PARATUCK-($N_1,N$) models in replacing $\mat A^{(1)}\otimes\mat A^{(2)}$, $\bm{\Phi}^{(1)}\diamond \bm{\Phi}^{(2)}$, $\mat C_{I_4\times R_1R_2}$, and $R_1R_2$, by $\overset{N_1}{\underset{n=1}{\otimes}} \mat A^{(n)}$, $\overset{N_1}{\underset{n=1}{\diamond}} \bm{\Phi}^{(n)}$, $\mat C_{I_{N_1+2}...I_N\times R}$, and $R=\overset{N_1}{\underset{n=1}{\prod}} R_n$, respectively.

\section{Conclusion}\label{sec:concl}

Several tensor models among which some are new, have been presented in a general and unified framework. {\color{black}The use of the index notation for mode combination based on Kronecker products provides an original and concise way to derive vectorized and matricized forms of tensor models.} A particular focus on constrained tensor models has been made with a perspective of designing MIMO communication systems with resource allocation. A link between PARATUCK models and constrained PARAFAC models has been established, which allows to apply results concerning PARAFAC models to derive uniqueness properties and parameter estimation algorithms for PARATUCK models. In a companion paper, several tensor-based MIMO systems are presented in a unified way based on constrained PARAFAC models, and a new tensor-based space-time-frequency (TSTF) MIMO transmission system with a blind receiver is proposed using a generalized PARATUCK model {\color{black}\cite{TSTF:TSP14}}. Even if this presentation of constrained tensor models has been made with the aim of designing MIMO transmission systems, we believe that such tensor models can be applied to other areas than telecommunications, like for instance biomedical signal processing, and more particularly for ECG and EEG signals modeling, with spatial constraints allowing to take into account the relative weight of the contributions of different areas of surface to electrodes. {\color{black} The considered constrained tensor models allow to take constraints into account either independently on each matrix factor of a PARAFAC decomposition, in the case of PARALIND/CONFAC models, or  between factors, in the case of PARATUCK models. A perspective of this work is to consider constraints into tensor networks which decompose high order tensors into lower-order tensors for big data processing \cite{Ci2014}. In this case, the constraints could act either separately on each tensor component to facilitate their physical interpretability, or between tensor components to explicit their interactions. }

\section*{Appendix}\label{sec:app}


\noindent {\bf{A1. Some matrix formulae}}\label{sec:appA}\\
For $\mat A^{(n)}\inmat{I_{n}}{R_{n}}$, $\mat B^{(n)}\inmat{R_{n}}{J_{n}}$, $\mat \Phi^{(n)}\inmat{R_{n}}{R}$, and $\mat \Psi^{(n)}\inmat{R_{n}}{Q}$, $n=1,\cdots,N$
\beq \label {eq: kron prop1}
&&\left( \overset{N}{\underset{n=1}{\kron}}\mat A^{(n)} \right)^T=\overset{N}{\underset{n=1}{\kron}}{\mat A^{(n)}}^T\inmat {R_1\cdots R_N}{I_1\cdots I_N}\\
\label {eq: kron assoc prop}
&&\left( \overset{N}{\underset{n=1}{\kron}}\mat A^{(n)} \right)\left( \overset{N}{\underset{n=1}{\kron}}\mat B^{(n)} \right)= \overset{N}{\underset{n=1}{\kron}}\mat A^{(n)}\mat B^{(n)}\inmat {I_1\cdots I_N}{J_1\cdots J_N}\\
\label {eq: kron prop3}
&&\left( \overset{N}{\underset{n=1}{\kron}}\mat A^{(n)} \right)\left( \overset{N}{\underset{n=1}{\diamond}}\mat \Phi^{(n)} \right)=\overset{N}{\underset{n=1}{\diamond}}\mat A^{(n)}\mat \Phi^{(n)}\inmat {I_1\cdots I_N}{R}\\
\label {eq: kron prop4}
&&\left( \overset{N}{\underset{n=1}{\diamond}}\mat \Psi^{(n)} \right)^T\left( \overset{N}{\underset{n=1}{\diamond}}\mat \Phi^{(n)} \right)=\overset{N}{\underset{n=1}{\odot}}{\mat \Psi^{(n)}}^T\mat \Phi^{(n)}\inmat{Q}{R}.\nonumber\\
&&\qquad \qquad \qquad (\textrm{Associative Property})\nonumber
\eeq
For $\mat A^{(n)}\inmat{I}{J}$, $n=1,\cdots,N$, and $\mat B^{(p)}\inmat{K}{L}$, $p=1,\cdots,P$
\beq \label {eq: kron distribut prop}
&&\hspace{-4ex}\left(\sum\limits_{n=1}^{N}\mat A^{(n)} \right)\kron\left(\sum\limits_{p=1}^{P}\mat B^{(p)} \right)=\sum\limits_{n=1}^{N}\sum\limits_{p=1}^{P}\left(\mat A^{(n)}\kron\mat B^{(p)} \right)\inmat{IK}{JL}\\ 
&&\qquad \qquad \qquad (\textrm{Distributive Property})\nonumber
\eeq
In particular, for $\mat A\inmat{I}{M}$, $\mat B\inmat{J}{N}$, $\mat C\inmat{M}{P}$, $\mat D\inmat{N}{Q}$, $\mat E\inmat{P}{J}$, $\mat \Phi\inmat{M}{R}$, $\mat \Psi\inmat{N}{R}$, $\mat \Omega\inmat{M}{Q}$, $\mat \Xi\inmat{N}{Q}$, and $\mat x\inmat{M}{1}$, we have
\beq \label{eq: kron prop5}
&&(\mat A\kron\mat B)^T={\mat A}^T\kron{\mat B}^T,\\
\label {eq: kron prop6}
&&(\mat A\kron\mat B)(\mat C\kron\mat D)={\mat A\mat C}\kron{\mat B\mat D},\\
\label{eq:Kron-Khatri}
&&(\mat A\kron\mat B)(\mat\Phi\diamond\mat\Psi)={\mat A\mat \Phi}\diamond{\mat B\mat \Psi},\\
&&(\mat\Omega\diamond\mat\Xi)^T(\mat\Phi\diamond\mat\Psi)={\mat \Omega}^T\mat\Phi\odot{\mat \Xi}^T\mat \Psi,\\
\label{eq:vec formula a}
&&\textrm{vec}(\mat A\mat C\mat E)=({\mat E}^T\kron\mat A)\textrm{vec}(\mat C),\\
\label{eq:vec formula b}
&&\textrm{vec}\left(\mat A diag(\mat x)\mat C \right)=({\mat C}^T\diamond \mat A)\mat x.
\eeq

\

\noindent {\bf{A2. Proof of (\ref{eq:Tucker mat rep})}}\\
{\color{black}
Defining ($\mathds{I}_1,\mathds{I}_2$) and ($\mathds{R}_1,\mathds{R}_2$) as the sets of indices $i_n$ and $r_n$ associated respectively with the sets ($\mathds{S}_1,\mathds{S}_2$) of index $n$, the formula (\ref{unfolinverse ind not}) allows writing the element $g_{r_1,\cdots,r_N }$ of the core tensor as
\beq \label {eq: core tensor element b}
g_{r_1,\cdots,r_N }=\mat e^{\mathds{R}_1}\mat G_{\mathds{S}_1;\mathds{S}_2}\mat e_{\mathds{R}_2}.
\eeq
where $\mathds{R}_1=\{r_n, n\in \mathds{S}_1\}$ and $\mathds{R}_2=\{r_n, n\in \mathds{S}_2\}$.\\
Substituting $x_{i_1,\cdots,i_N }$ and $g_{r_1,\cdots,r_N }$ by their expressions (\ref{Tucker ind not}) and (\ref {eq: core tensor element b}) into (\ref {unfol ind not}) gives
\beq
\mat X_{\mathds{S}_1;\mathds{S}_2}&=&x_{i_1,\cdots,i_N }\mat e_{\mathds{I}_1}^{\mathds{I}_2}=\mat e_{\mathds{I}_1}x_{i_1,\cdots,i_N }\mat e^{\mathds{I}_2}\nonumber\\
&=&\mat e_{\mathds{I}_1}g_{r_1,\cdots,r_N }\prod\limits_{n=1}^{N}a^{(n)}_{i_n,r_n}\mat e^{\mathds{I}_2}\nonumber\\
&=&\prod\limits_{n=1}^{N}a^{(n)}_{i_n,r_n}\mat e_{\mathds{I}_1}\mat e^{\mathds{R}_1}\mat G_{\mathds{S}_1;\mathds{S}_2}\mat e_{\mathds{R}_2}\mat e^{\mathds{I}_2}\nonumber\\
&=&(\underset{n\in \mathds{S}_1}{\prod}a^{(n)}_{i_n,r_n}\mat e_{\mathds{I}_1}^{\mathds{R}_1})\mat G_{\mathds{S}_1;\mathds{S}_2}(\underset{n\in \mathds{S}_2}{\prod}a^{(n)}_{i_n,r_n}\mat e_{\mathds{R}_2}^{\mathds{I}_2})
\eeq
Applying the general Kronecker formula (\ref{Kron ind not}) in terms of the index notation allows to rewrite this matrix unfolding as
\beq
\mat X_{\mathds{S}_1;\mathds{S}_2}=\left( \underset{n \in \mathds{S}_1}{\otimes} \mat A^{(n)} \right)\mat G_{\mathds{S}_1;\mathds{S}_2}\left( \underset{n \in \mathds{S}_2}{\otimes} \mat A^{(n)} \right)^T. \nonumber
\eeq
}

\

\noindent {\bf{A3. Proof of (\ref{eq:Parafac mat rep})}}\\
{\color{black}
Substituting the expression (\ref {PARAFAC ind not}) of $x_{i_1,\cdots,i_N}$ into (\ref {unfol ind not}) and using the identities (\ref{veckhatri ind not}) and (\ref{UV ind not}) give
\beq
\mat X_{\mathds{S}_1;\mathds{S}_2}&=&x_{i_1,\cdots,i_N}\mat e_{\mathds{I}_1}^{\mathds{I}_2}\nonumber\\
&=&\left(\underset{n\in \mathds{S}_1}{\prod}a_{i_n,r}^{(n)}\mat e_{\mathds{I}_1} \right)\left(\underset{n\in \mathds{S}_2}{\prod}a_{i_n,r}^{(n)}\mat e^{\mathds{I}_2} \right)\nonumber\\
&=&\left( \underset{n \in \mathds{S}_1}{\kron}\mat A^{(n)}_{.r} \right)\left( \underset{n \in \mathds{S}_2}{\kron}\mat A^{(n)}_{.r} \right)^T\nonumber\\
&=&\left( \underset{n \in \mathds{S}_1}{\diamond}\mat A^{(n)} \right)\left( \underset{n \in \mathds{S}_2}{\diamond}\mat A^{(n)} \right)^T
\eeq
which ends the proof of (\ref{eq:Parafac mat rep}).
}

\

\noindent {\bf{A4. Proof of (\ref{eq:Const Parafac factors}) and (\ref{eq:Paratuck-(2,4) constraints})}}\label{sec:appB}\\
Let us define the third-order tensors $\matc A \inten{I_1}{R_1}{R_2}$, $\matc B \inten{I_2}{R_1}{R_2}$, $\matc F\inten{I_3}{R_1}{R_2}$, and $\matc D \inten{I_4}{R_1}{R_2}$ such as
\beq
a_{i_1,r_1,r_2}=a^{(1)}_{i_1,r_1}\,\, \forall r_2=1,\cdots, R_2 \,&;&\quad b_{i_2,r_1,r_2}=a^{(2)}_{i_2,r_2}\,\, \forall r_1=1,\cdots, R_1;\nonumber\\
\label{eq:var def}
f_{i_3,r_1,r_2}=\phi^{(1)}_{r_1,i_3}\phi^{(2)}_{r_2,i_3} \,&;&\quad d_{i_4,r_1,r_2}=c_{r_1,r_2,i_4}.
\eeq
\
The tensor model (\ref{eq:Paratuck-(2,4) model}) can be rewritten as
\beq \label{Paratuck-(2,4) model equiv1}
x_{i_1,i_2,i_3,i_4}=\sum\limits_{r_1=1}^{R_1}\sum\limits_{r_2=1}^{R_2}a_{i_1,r_1,r_2}b_{i_2,r_1,r_2}f_{i_3,r_1,r_2}d_{i_4,r_1,r_2}.
\eeq
\
Defining the change of variables $r=(r_1-1)R_2+r_2$ that corresponds to a combination of the last two modes of the tensors $\matc A$, $\matc B$, $\matc F$, and $\matc D$, {\color{black} Eq.
 (\ref{Paratuck-(2,4) model equiv1}) can be rewritten as the constrained PARAFAC-4 model (\ref{eq:Paratuck-(2,4) model equiv}), where $a_{i_1,r}$, $b_{i_2,r}$, $f_{i_3,r}$, and $d_{i_4,r}$ are entries of mode-1 matrix unfoldings of the tensors $\matc A$, $\matc B$, $\matc F$, and $\matc D$, i.e. entries of $\mat A \overset{\vartriangle}{=}\mat A_{I_1\times R_1R_2}$, $\mat B\overset{\vartriangle}{=}\mat B_{I_2\times R_1R_2}$, $\mat F\overset{\vartriangle}{=}\mat F_{I_3\times R_1R_2}$, and $\mat D\overset{\vartriangle}{=}\mat D_{I_4\times R_1R_2}$, respectively.
 Using the formulae (\ref{tensor extension 1}) and (\ref{tensor extension 2}), we can directly deduce the following expressions of $\mat A$ and $\mat B$
\beq
\mat A&=&\mat A^{(1)}\otimes \mat 1^T_{R_2}=\mat A^{(1)}(\mat I_{R_1}\otimes {\mat 1}^T_{R_2})=\mat A^{(1)}\mat \Psi^{(1)}.\label{eq:formula APsi1}\\
\mat B&=&\mat 1^T_{R_1}\otimes \mat A^{(2)}=\mat A^{(2)}(\mat 1^T_{R_1}\otimes \mat I_{R_2})=\mat A^{(2)}\mat \Psi^{(2)}\label{eq:formula APsi2}
\eeq
For the matrix ${\mat F}$, using the index notation with the definition (\ref{eq:var def}) gives
\beq
\mat F=(f_{i_3,r_1,r_2}\mat e^{r_1r_2}_{i_3})=(\phi^{(1)}_{r_1,i_3} \phi^{(2)}_{r_2,i_3}\mat e^{r_1r_2}_{i_3})\nonumber
\eeq
Applying the formula (\ref{Khatri index notation}), we directly obtain
\beq
\mat F=(\mat \Phi^{(1)}\diamond \mat \Phi^{(2)})^T.\nonumber
\eeq
}

\

{\color{black}
\noindent {\bf{A5. Tensor extension of a matrix}}}
Following the same demonstration as for (\ref{tensor extension 1}) and (\ref{tensor extension 2}), it is easy to deduce the following more general formula for the extension of $\mat B\inmat{I}{R_n}$ into a tensor $\matc A\inten{I}{R_1\times \cdots}{R_N}$ such as $a_{i,r_1,\cdots,r_n,\cdots,r_N}=b_{i,r_n} \,\, \forall \, r_k=1,\cdots,R_k, \,\, \textrm{for} \,\, k=1,\cdots,n-1,n+1,\cdots,N$. Defining $R=\prod\limits_{n=1}^{N}R_n$, we have
\beq \label{eq: extension formula}
\hspace{-5ex}\mat A_{I\times R}&=&\mat B ({\mat 1}^T_{R_1}\otimes\cdots\otimes {\mat 1}^T_{R_{n-1}}\otimes\mat I_{R_n}\otimes{\mat 1}^T_{R_{n+1}}\otimes\cdots \otimes{\mat 1}^T_{R_N})\inmat{I}{R}.
\eeq
{\color{black}
Similarly, for the extension of $\mat B\inmat{I_n}{R}$ into a tensor $\matc A\inten{I_1\times \cdots}{I_N}{R}$ such as $a_{i_1,\cdots,i_n,\cdots,i_N,r}=b_{i_n,r} \,\, \forall \, i_k=1,\cdots,I_k, \,\, \textrm{for} \,\, k=1,\cdots,n-1,n+1,\cdots,N$, we have
\beq \label{eq: extension formula b}
\hspace{-5ex}\mat A_{I\times R}&=& ({\mat 1}_{I_1}\otimes\cdots\otimes {\mat 1}_{I_{n-1}}\otimes\mat I_{I_n}\otimes{\mat 1}_{I_{n+1}}\otimes\cdots \otimes{\mat 1}_{I_N})\mat B\inmat{I}{R}.
\eeq
where $I=\prod\limits_{n=1}^{N}I_n$.\\
For instance, if we consider the following tensor extension of $\mat B\inmat{I}{J}$
\be
a_{m,n,i,j,k,l}=b_{i,j} \,\, \forall \, m=1,\cdots,M, \forall \, n=1,\cdots,N, \forall \, k=1,\cdots,K, \forall \, l=1,\cdots,L \nonumber
\ee
\
the combination of formulae (\ref{eq: extension formula}) and (\ref{eq: extension formula b}) gives
\be \label{eq: extension formula c}
\mat A_{MNI\times JKL}=(\mat 1_{MN}\otimes\mat I_I)\mat B (\mat I_J\otimes\mat 1^T_{KL})
\ee
which can be written as
\be
\mat A_{MNI\times JKL}=\mat B \times_1 \bm{\Psi}_1 \times_2 (\bm{\Psi}_2)^T\nonumber
\ee
with $\bm{\Psi}_1=\mat 1_{MN}\otimes\mat I_I$ and $\bm{\Psi}_2=\mat I_J\otimes\mat 1^T_{KL}$.
}

\section*{Acknowledgements}
This work has been developed under the FUNCAP/CNRS bilateral cooperation project (2013-2014).\\
Andr\'{e} L. F. de Almeida is partially supported by CNPq. The authors are thankful to A. Cichocki for  useful comments and suggestions.


\bibliographystyle{bmc-mathphys} 

\begin{thebibliography}{10}

\bibitem{Cullagh}
P McCullagh, \emph{Tensor methods in statistics}. (Chapman and Hall, London, New York, 1987)

\bibitem{Como02:oxford}
P Comon, Tensor decompositions: State of the art and applications, in
\emph{Mathematics in Signal Processing {V}}, JG McWhirter and IK
  Proudler, Eds. Oxford, UK: Clarendon Press, 1--24 (2002)

\bibitem{Tucker66}
LR Tucker, Some mathematical notes on three-mode factor analysis,
Psychometrika, \textbf{31}, 279--311 (1966)

\bibitem{Harshman70}
RA Harshman, Foundations of the {PARAFAC} procedure: Model and conditions
  for an ``explanatory'' multimodal factor analysis, {UCLA} Working
  Papers in Phonetics, \textbf{16}, 1--84, (1970)

\bibitem{CarrollChang70}
JD Carroll and J Chang, Analysis of individual differences in
  multidimensional scaling via an {N}-way generalization of ``{E}ckart-{Y}oung''
  decomposition, Psychometrika, \textbf{35}(3), 283--319 (1970)

\bibitem{Kiers2000}
HAL Kiers, Towards a standardized notation and terminology in multiway
  analysis, J. Chemometrics, \textbf{14}(2), 105--122 (2000)

\bibitem{Kroonenberg2008}
PM Kroonenberg, \emph{Applied multiway data analysis}. (John Wiley and Sons, 2008)

\bibitem{Bro98}
R Bro, \emph{Multi-way analysis in the food industry: {M}odels, algorithms and
  applications}, Ph.D. dissertation, University of Amsterdam, Amsterdam (1998)

\bibitem{Smilde_book}
A Smilde, R Bro, and P Geladi, \emph{Multi-way Analysis. Applications in the
  Chemical Sciences}. (John Wiley and Sons, Chichester, UK, 2004)

\bibitem{Cardoso90}
J-F Cardoso, in  IEEE ICASSP'90. Eigen-structure of the fourth-order cumulant tensor with
  application to the blind source separation problem (Albuquerque, USA, 1990), pp. 2655--2658

\bibitem{CardoComon90}
J-F Cardoso and P Comon, in EUSIPCO'90. Tensor-based independent component analysis
(Barcelona, Spain, 1990), pp. 673--676.

\bibitem{Cardoso91}
J-F Cardoso, in Proceedings of IEEE ICASSP'91, Super-symmetric decomposition of the fourth-order cumulant
  tensor. {B}lind identification of more sources than sensors (Toronto, Canada, 1991), pp. 3109--3112

\bibitem{DeLathau97}
L De Lathauwer, \emph{Signal processing based on multilinear algebra}, Ph.D.
  dissertation, KU Leuven, Leuven, 1997.

\bibitem{CJ10}
P Comon and C Jutten, \emph{Handbook of blind source separation. Independent
  component analysis and applications} (Elsevier, Oxford, UK, 2010)

\bibitem{Sid00}
ND Sidiropoulos, GB Giannakis, and R Bro, Blind {PARAFAC} receivers
  for {DS-CDMA} systems, {IEEE} Trans. Signal Process., \textbf{48}(3), 810--823 (2000)

\bibitem{Cichocki09}
A Cichocki, R Zdunek, AH Phan, and S-I Amari, \emph{Nonnegative matrix
  and tensor factorizations. Applications to exploratory multi-way data
  analysis and blind source separation}, (Wiley, Chichester, UK, 2009)

\bibitem{Kolda09}
TG Kolda and BW Bader, Tensor decompositions and applications, {SIAM} J. Matrix Anal. Appl., \textbf{51}(3), 455--500 (2009)

\bibitem{Acar09}
E Acar and B Yener,  Unsupervised multiway data analysis: {A} literature
  survey, IEEE Trans. Knowledge and data engineering, \textbf{21}(1), 6--20 (2009)

\bibitem{Morten11}
M Morup, Applications of tensor (multiway array) factorizations and decompositions in data mining,
Wiley Interdisciplinary Reviews: Data mining and knowledge discovery, John Wiley and Sons,
\textbf{1}(1), 24--40 (2011).

{\color{black}
\bibitem{SLC2012}
M Sorensen, L De Lathauwer, P Comon, S Icart, and L Deneire, Canonical polyadic decomposition with a columnwise orthonormal factor matrix, SIAM J. Matrix Analysis and Appl., \textbf{33}(4), 1190-1213 (2012).
}

\bibitem{Shashua2005}
A Shashua and T Hazan, in Proc. of 22nd Int. Conf. on
  {M}achine {L}earning, Non-negative tensor factorization with applications
  to statistics and computer vision, (Bonn, Germany, 2005), pp. 792--799

\bibitem{Hazan2005}
S Hazan, S Polak, and A Shashua, in Proc. of 10th IEEE Int. Conf.
  on Computer Vision ({ICCV'2005}), Sparse image coding using a 3{D}
  non-negative tensor factorization, (Beijing, China, 2005), pp. 50--57

\bibitem{Friedlander2008}
MP Friedlander and K Hatz, Computing nonnegative tensor factorizations,
  Optimization Methods and Software, \textbf{23}(4), 631--647 (2008)

\bibitem{Zhang2008}
Q Zhang, H Wang, R Plemmons, and P Pauca, Tensor methods for
  hyperspectral data processing: {A} space object identification study,
  J. Opt. Soc. Am. A, \textbf{25}(12), 3001--3012 (2008)

\bibitem{Benetos2010}
E Benetos and C Kotropoulos, Non-negative tensor factorization applied to
  music genre classification, {IEEE} Trans. on Audio, Speech, and
  Language Proc., \textbf{18}(8), 1955--1967 (2010)

\bibitem{Ozerov2011}
A Ozerov, C F\'{e}vote, R Blouet, and G Durrieu, in International Conference on Acoustics, Speech and
  Signal Processing (ICASSP2011), Multichannel nonnegative
  tensor factorization with structured constraints for user-guided audio source
  separation, (Prague, Czech Republic, 2011)

\bibitem{Acar2010}
E Acar, DM Dunlavy, TG Kolda, and M Morup, in Proc. of 10th {SIAM} Int. Conf.
  on Data mining, Scalable tensor factorizations with missing data,
  (Columbus, Ohio, 2010), pp. 701--712

\bibitem{Royer2011}
J-P Royer, N Thirion-Moreau, and P Comon, Computing the polyadic
  decomposition of nonnegative third order tensors, Signal Processing,
  \textbf{91}, 2159--2171 (2011)

\bibitem{Phan2011}
A-H. Phan and A Cichocki, Extended {HALS} algorithm for nonnegative
  {T}ucker decomposition and its applications for multiway analysis and
  classification, Neurocomputing, \textbf{74}, 1956--1969 (2011)

\bibitem{Welling2001}
M Welling and M Weber, Positive tensor factorization, Pattern
  Recogn. Letters, \textbf{22}(12), 1255--1261 (2001)

\bibitem{Morup2008}
M Morup and LK Hansen, Algorithms for sparse non-negative {T}ucker
  decompositions, Neural Computation, \textbf{20}, 2112--2131 (2008)

\bibitem{FB09EUSIPCO}
G Favier and T Bouilloc, in European Sign. Proc. Conf. (EUSIPCO'2010),
A constrained tensor based approach for {MIMO NL-CDMA} systems,
(Aalborg, Denmark, 2010)

\bibitem{FBA12}
G Favier, T Bouilloc, and ALF de~Almeida, Blind constrained
  block-{T}ucker2 receiver for multiuser {SIMO NL-CDMA} communication
  systems, Signal Processing, \textbf{92}(7), 1624--1636 (2012)

\bibitem{FKB2012}
G Favier, AY Kibangou, and T Bouilloc, Nonlinear system modeling and
  identification using {V}olterra-{PARAFAC} models, Int. J. of Adaptive
  Control and Sig. Proc., \textbf{26}, 30--53 (2012)

\bibitem{BF2012}
T Bouilloc and G Favier,  Nonlinear channel modeling and identification
  using bandpass {V}olterra-{PARAFAC} models, Signal Processing,
 \textbf{92}(6), 1492--1498 (2012)

\bibitem{KF2009}
AY Kibangou and G Favier, Identification of parallel-cascade {W}iener
  systems using joint diagonalization of third-order {V}olterra kernel
  slices, {IEEE} Signal Proc. Letters, \textbf{16}(3) (2009)

\bibitem{Favier2009}
G Favier, in Proc. of 10th Int. Conf. on Sciences and Techniques of
  Automatic Control and Computer Engineering (STA'2009),
  Nonlinear system modeling and identification using tensor approaches,
  (Hammamet, Tunisia,2009)

\bibitem{Fernandes2007}
CER Fernandes, G Favier, and JCM Mota, Blind channel
  identification algorithms based on the {P}arafac decomposition of cumulant
  tensors: {T}he single and multiuser cases, Signal Processing,
  \textbf{88}, 1382--1401 (2008)

\bibitem{Fernandes2009}
CER Fernandes, G Favier, and JCM Mota, in Proc. of 15th
  {IFAC} Symp. on System Identification ({SYSID'2009}), Parafac-based blind
  identification of convolutive {MIMO} linear systems, (Saint-Malo, France, 2009)

\bibitem{Brachat2010}
J Brachat, P Comon, B Mourrain, and E Tsigaridas, Symmetric tensor
  decomposition, Linear Algebra and its Appl., \textbf{433}(11-12),
  1851--1872 (2010)

\bibitem{Nion2008}
D Nion and L {D}e Lathauwer, A block component model-based blind {DS-CDMA}
  receiver, {IEEE} Trans. Signal Proc., \textbf{56}(11), 5567--5579 (2008)

\bibitem{KF2009b}
AY Kibangou and G Favier, in European Signal Proc. Conf. (EUSIPCO'2009),
Noniterative solution for {P}arafac with a {T}oeplitz factor, (Glasgow, UK, 2009)
{\color{black}
\bibitem{SL2013}
M Sorensen, and L De Lathauwer, Blind signal separation via tensor decomposition with Vandermonde factor: canonical polyadic decomposition,
  IEEE Trans. Signal Process., \textbf{61}(22), 5507-5519 (2013).}
{\color{black}
\bibitem{GoF2014}
JH Goulart, and G Favier, An algebraic solution for the CANDECOMP/PARAFAC decomposition with circulant factors,
Submitted to Linear Algebra and its Applications (Feb. 2014)
http://hal.archives-ouvertes.fr/docs/00/96/72/63/PDF/RR-2014-02{\_}I3S.pdf.}

\bibitem{Comon2010}
P Comon, M Sorensen, and E Tsigaridas, in Proc. of IEEE
  ICASSP'2010, Decomposing tensors with structured matrix factors reduces to rank-1 approximations,
  (Dallas, USA, 2010), pp. 14--19.
{\color{black}
\bibitem{SC2013}
M Sorensen, and P Comon, Tensor decompositions with banded matrix factors,
Linear Algebra and its Applications, \textbf{438}, 919-941 (2013).}

\bibitem{CPK80}
JD Carroll, S Pruzansky, and JB Kruskal, Candelinc: a general approach
  to multidimensional analysis of many-way arrays with linear constraints on
  parameters, Psychometrika, \textbf{45}(1), 3--24 (1980).

{\color{black}
\bibitem{P2011}
DSG Pollock, On Kronecker products, tensor products and matrix differential calculus, Working paper 11/34,
Univ. of Leicester, Dept. of Economics, UK, {http://www.le.ac.uk/ec/research/RePEc/lec/leecon/dp11-34.pdf}, (July 2011).
}
%

\bibitem{DeLathau00}
L {D}e Lathauwer, B {D}e Moor, and  J {V}andewalle, A multilinear singular
  value decomposition, {SIAM} J. Matrix Anal. Appl., \textbf{21}(4), 1253--1278 (2000)

\bibitem{Cattell44}
RB Cattell, Parallel proportional profiles, and other principles for
  determining the choice of factors by rotation, Psychometrika,
 \textbf{9}, 267--283 (1944)

\bibitem{Hitchcock27}
FL Hitchcock, The expression of a tensor or a polyadic as a sum of
  products, Journal of Mathematics and Physics, 6(3), 164--189 (1927)

\bibitem{Kru77}
JB Kruskal, Three-way arrays: Rank and uniqueness of trilinear
  decompositions, with application to arithmetic complexity and statistics,
  Linear Algebra Appl., \textbf{18}(2), 95--138 (1977)

\bibitem{Com09:LinearAlgebra}
P Comon, JMF ten Berge, L De Lathauwer, and J Castaing, Generic and
  typical ranks of multi-way arrays, Linear Algebra and its
  Applications, \textbf{430}(11), 2997--3007 (2009)

\bibitem{ComoGLM08:SIAM}
P Comon, G Golub, L-H Lim, and B Mourrain, Symmetric tensors and
  symmetric tensor rank, SIAM J. Matrix Anal. Appl., \textbf{30}(3),
  1254--1279 (2008)

\bibitem{Bro2003}
R Bro and HAL Kiers, A new efficient method for determining the number
  of components in {PARAFAC} models, J. Chemometrics, \textbf{17}(5),
  274--286 (2003)

\bibitem{daCosta2008}
JPCL da Costa, M Haardt, and F Roemer, in Proc. of 5th {IEEE} Sensor Array and Multich. Signal Proc. Workshop (SAM 2008),
Robust methods based on {HOSVD} for estimating the model order in {PARAFAC} models, (Darmstadt, Germany, 2008), pp. 510--514

\bibitem{daCosta2010}
JPCL da Costa, F Roemer, M Weis, and M Haardt, in Proc. of {ITG}
  Workshop on Smart Antennas (WSA 2010), Robust ${R}$-{D}
  parameter estimation via closed-form {PARAFAC},  (Bremen, Germany, 2010), pp. 99--106

\bibitem{daCosta2011}
JPCL da Costa, F Roemer, M Haardt, and RT de Sousa, Multi-dimensional model order selection, EURASIP J. on Advances in
  Signal Processing, \textbf{26}, (July 2011)

\bibitem{Harshman96}
RA Harshman and ME Lundy, Uniqueness proof for a family of models
  sharing features of {T}ucker's three-mode factor analysis and
  {PARAFAC/CANDECOMP}, Psychometrika, \textbf{61}, 133--154 (1996)

\bibitem{Bro05}
R Bro, RA Harshman, and ND Sidiropoulos, Modeling multi-way data with
  linearly dependent loadings, {KVL} tech. report 176, (2005)

\bibitem{BHSL09}
R Bro, RA Harshman, ND Sidiropoulos, and ME Lundy, Modeling
  multi-way data with linearly dependent loadings, Chemometrics,
  \textbf{23}(7-8), 324--340 (2009)

\bibitem{Kibangou07Eusipco}
AY Kibangou and G Favier, in Proc. of European Signal Processing Conference (EUSIPCO'2007), Blind joint identification and equalization of
  {W}iener-{H}ammerstein communication channels using {PARATUCK-2} tensor decomposition,
  (Poznan, Poland, Sept. 2007)

\bibitem{Xu2012a}
L Xu, J Ting, Y Longxiang, and Z Hongbo, {PARALIND}-based identifiability
  results for parameter estimation via uniform linear array, EURASIP J.
  Advances in Sig. Proc. (2012)

\bibitem{Xu2012b}
L Xu, G Liang, Y Longxiang, and Z Hongbo, {PARALIND}-based blind joint
  angle and delay estimation for multipath signals with uniform linear array,
  EURASIP J. Advances in Sig. Proc. (2012)

\bibitem{AFM08IEEESP}
ALF de Almeida, G Favier, and JCM Mota, A constrained factor
  decomposition with application to {MIMO} antenna systems, IEEE Trans.
  Signal Process., \textbf{56}(6), 2429--2442 (2008)

\bibitem{FCAR11EUSIPCO}
G Favier, MN da~Costa, ALF de~Almeida, and JMT Romano, in Proc. of European
  Sign. Proc. Conf. (EUSIPCO'2011), Tensor
 coding for {CDMA-MIMO} wireless communication systems, (Barcelona, Spain, Aug. 29-Sept. 2 2011)

\bibitem{FCAR12SP}
G Favier, MN da~Costa, ALF de~Almeida, and JMT Romano, Tensor
  space-time ({TST}) coding for {MIMO} wireless communication systems,
  Signal Processing, \textbf{92}(4), 1079--1092 (2012)

\bibitem{Almeida2009}
ALF de~Almeida, G Favier, and JCM Mota, Space-time
  spreading-multiplexing for {MIMO} wireless communication systems using the
  {PARATUCK-2} tensor model, Signal Processing, \textbf{89}(11), 2103--2116 (Nov. 2009)

\bibitem{ALSC2012}
ALF de~Almeida, X Luciani, A Stegeman, and P Comon, CONFAC decomposition approach to blind
identification of underdetermined mixtures based on generating function derivatives,
IEEE Trans. Signal Process., \textbf{60}(11), 5698-5713 (2012).

\bibitem{Almeida05Asilomar}
ALF de~Almeida, G Favier, and JCM Mota, in
  Asilomar Conf. Sig. Syst. Comp., Generalized {PARAFAC} model for multidimensional wireless
  communications with application to blind multiuser equalization, (Pacific Grove, CA, USA, Nov. 2005)

\bibitem{Almeida05PSIP}
ALF de~Almeida, G Favier, and JCM Mota, in Int.  Conf. on Physics in Signal and Image processing (PSIP), {PARAFAC} models for wireless communication systems, (Toulouse, France, Jan. 31 - Feb. 2, 2005)

\bibitem{Almeida07SP}
ALF de~Almeida, G Favier, and JCM Mota, {PARAFAC}-based unified tensor modeling for wireless communication systems with application to blind multiuser equalization, Signal
 Processing, \textbf{87}, 337--351 (2007).

\bibitem{Almeida06PIMRC}
ALF de~Almeida, G Favier, and JCM Mota, in Proc. of 17th {IEEE} Symp. Pers. Ind. Mob.
  Radio Com. (PIMRC'2006), Tensor-based space-time multiplexing codes for {MIMO-OFDM} systems
  with blind detection, (Helsinki, Finland, Sept. 2006).

\bibitem{Almeida2009b}
ALF de~Almeida, G Favier, and JCM Mota, Constrained {T}ucker-3 model for blind beamforming,
Signal  Processing, \textbf{89}, 1240--1244 (2009).

{\color{black}
\bibitem{SRK2009}
J Salmi, A Richter, and V Koivunen, Sequential unfolding SVD for tensors with applications in array signal processing, IEEE Trans. Signal Process. \textbf{57}(12), 4719-4733 (Dec. 2009).
}

\bibitem{DeLathauwer2008}
L {D}e Lathauwer, Decompositions of a higher-order tensor in block
  terms-{P}art {II}: {D}efinitions and uniqueness, {SIAM} J. Matrix
  Anal. Appl., \textbf{30}(3), 1033--1066 (2008).

{\color{black}
\bibitem{CMP2014}
A Cichocki, D Mandic, A-H Phan, C Caiafa, G Zhou, Q Zhao, and L {D}e Lathauwer, Tensor decompositions for signal processing applications. From two-way to multiway component analysis, IEEE Signal Processing Magazine (to appear), arXiv:1403.4462v1 (March 2014)
}

\bibitem{AFM08IEEESP:2}
ALF de~Almeida, G Favier, and JCM Mota, Constrained tensor modeling
  approach to blind multiple-antenna {CDMA} schemes, IEEE Trans. Signal
  Process., \textbf{56}(6), 2417--2428 (June 2008).

{\color{black}
\bibitem{L2011}
L de~Lathauwer, Blind separation of exponential polynomials and the decomposition of a tensor in rank-$(L_r,L_r,1)$ terms, SIAM J. Matrix Anal. Appl., \textbf{32}(4), 1451-1474 (2011).
}

\bibitem{Sid00Nway}
ND Sidiropoulos, and R Bro, On the uniqueness of multilinear decomposition of {N}-way arrays, J. Chemometrics, \textbf{14}, 229--239 (2000)

{\color{black}
\bibitem{DL2013}
I Domanov, and L De Lathauwer, On the uniqueness of the canonical polyadic decomposition of third-order tensors. Part I: Basic results and uniqueness of one factor matrix, arXiv:1301.4602v1, KU Leuven, Belgium, (Jan 2013).

\bibitem{DL2014}
I Domanov, and L De Lathauwer, Generic uniqueness conditions for the canonical polyadic decomposition and INDSCAL, arXiv:1405.6238v1, KU Leuven, Belgium, (May 2014).
}

\bibitem{tenBerge022}
JMF ten Berge, and ND Sidiropoulos, On uniqueness in {CANDECOMP/PARAFAC}, Psychometrika, \textbf{67}(3), 399-409 (2002)

\bibitem{Harshman72}
RA Harshman, Determination and proof of minimum uniqueness conditions for {PARAFAC1}, {UCLA} Working Papers in Phonetics, \textbf{22}, 111-117 (1972)

\bibitem{Stegeman07}
A Stegeman, and ND Sidiropoulos, On {K}ruskal's uniqueness condition for the {CANDECOMP/PARAFAC} decomposition, Lin. Alg. Appl., \textbf{420}, 540--552 (2007)

\bibitem{JiangSid04}
T Jiang, and ND Sidiropoulos, {K}ruskal's permutation lemma and the identification of {CANDECOMP/PARAFAC} and bilinear models with constant modulus constraints, {IEEE} Trans. Signal Process., \textbf{52}(9), 2625-2636 (Sept. 2004)

\bibitem{DeLathauwer2006}
L {De} Lathauwer, A link between the canonical decomposition in multilinear algebra and simultaneous matrix diagonalization, SIAM J. Matrix Anal. Appl., \textbf{28}(3), 642-666 (2006)

\bibitem{Stegeman08}
A Stegeman, On uniqueness conditions for {CANDECOMP/PARAFAC} and {INDSCAL} with full column rank in one mode, Lin. Alg. Appl., \textbf{431}(1-2), 211-227 (2008)

\bibitem{GMB11}
X Guo, S Miron, D Brie, S Zhu, and X Liao, A {CANDECOMP/PARAFAC} perspective on uniqueness of {DOA} estimation using a vector sensor array, {IEEE} Trans. Signal Process.,  \textbf{59}(7), 3475-3481 (July 2011)

{\color{black}

\bibitem{B2004}
JMF Ten Berge, Partial uniqueness in CANDECOMP/PARAFAC, J. Chemometrics, \textbf{18}, 12-16 (2004).

\bibitem{BMC2011}
D Brie, S Miron, F Caland, and C Mustin, in IEEE ICASSP'2011, An uniqueness condition for the 4-way CANDECOMP/PARAFAC model with collinear loadings in three modes, (Prague, Czech Republic, May 2011).

\bibitem{StegAlm09}
A. Stegeman, ALF de Almeida, Uniqueness conditions for constrained three-way factor decompositions with linearly dependent loadings, SIAM. J. Matrix Anal. Appl., \textbf{31}(3), 1469--1490 (Dec 2009).

\bibitem{StegLam12}
A Stegeman, and TTT Lam, Improved uniqueness conditions for canonical tensor decompositions with linearly dependent loadings, SIAM. J. Matrix Anal. Appl., \textbf{33}(4), 1250--1271 (Nov 2012).

\bibitem{GuoMironBrieSteg12}
X Guo, S Miron, D Brie, A Stegeman, Uni-mode and partial uniqueness conditions for CANDECOMP/PARAFAC of three-way arrays with linearly dependent loadings, SIAM. J. Matrix Anal. Appl., \textbf{33}(1), 111--129 (Jan 2012).

\bibitem{VanLoanPitsianis}
CF Van Loan, N Pitsianis, in \emph{Linear Algebra for Large Scale and Real-Time Applications} ed. MS Moonen, GH Golub, BLR de Moor (Kluwer Publications, Netherlands, 1993), p. 293.

\bibitem{TSTF:TSP14}
G Favier, and ALF de~Almeida,  Tensor space-time-frequency coding with semi-blind receivers for MIMO wireless communication systems, submitted to {IEEE} Trans. Signal Process, (Feb. 2014).

}

{\color{black}
\bibitem{Ci2014}
A Cichocki, Era of big data processing: A new approach via tensor networks and tensor decompositions, arXiv:1403.2048v3 (June 2014).
}

\end{thebibliography}


\end{document}